\documentclass[sn-basic]{sn-jnl}

\usepackage{graphicx}%
\usepackage{multirow}%
\usepackage{amsmath,amssymb,amsfonts}%
\usepackage{amsthm}%
\usepackage{mathrsfs}%
\usepackage[title]{appendix}%
\usepackage{xcolor}%
\usepackage{textcomp}%
\usepackage{manyfoot}%
\usepackage{booktabs}%
\usepackage{algorithm}%
\usepackage{algorithmic}
\usepackage{listings}%
\usepackage{tikz}
\usetikzlibrary{patterns}
\usepackage{pgfplots}
\usepgfplotslibrary{fillbetween}
\pgfplotsset{compat=1.13}

\theoremstyle{plain}
\newtheorem{theorem}{Theorem}[section]
\newtheorem{proposition}[theorem]{Proposition}
\newtheorem{lemma}[theorem]{Lemma}
\newtheorem{corollary}[theorem]{Corollary}
\theoremstyle{definition}
\newtheorem{definition}[theorem]{Definition}

\theoremstyle{remark}

\def\eqd{\stackrel{\mbox{\scriptsize{d}}}{=}}
\def\simiid{\stackrel{\mbox{\scriptsize{\rm iid}}}{\sim}}

\newcommand{\ddr}{\mathrm{d}}
\newcommand{\E}{\mathbb{E}}
\newcommand{\W}{\mathcal{W}}
\newcommand{\R}{\mathbb{R}}
\newcommand{\X}{\mathbb{X}}
\newcommand{\Y}{\mathbb{Y}}
\newcommand{\law}{\mathcal{L}}

\newcommand{\dlip}{d_{\rm{Lip}}}

\newcommand{\ngr}{M}
\renewcommand{\algorithmiccomment}[1]{\hfill// #1}

\begin{document}

\title[HIPM]{\center{Hierarchical Integral Probability Metrics:\\ A distance on random probability measures with \\low sample complexity}}


\author[1]{\fnm{Marta} \sur{Catalano}}

\author[2]{\fnm{Hugo} \sur{Lavenant}}

\affil[1]{\orgname{Luiss University}, \country{Italy}}

\affil[2]{\orgname{Bocconi University}, \country{Italy}}

\abstract{Random probabilities are a key component to many nonparametric methods in Statistics and Machine Learning. To quantify comparisons between different laws of random probabilities several works are starting to use the elegant Wasserstein over Wasserstein distance. In this paper we prove that the infinite dimensionality of the space of probabilities  drastically deteriorates its sample complexity, which is slower than any polynomial rate in the sample size. We propose a new distance that preserves many desirable properties of the former while achieving a parametric rate of convergence. In particular, our distance 1) metrizes weak convergence; 2) can be estimated numerically through samples with low complexity; 3) can be bounded analytically from above and below. The main ingredient are integral probability metrics, which lead to the name \emph{hierarchical IPM}.}


\maketitle

\section{Introduction}
\label{sec:intro}

In this work we discuss distances between laws of random probabilities, that is, probability measures on spaces of probability measures. Our motivation comes from nonparametric methods in Statistics and Machine Learning, often referred to as distribution-free in that they do not rely on strong assumptions on the underlying distribution of the data.
The need for flexibility has inspired a wealth of both frequentist and Bayesian nonparametric models, which often share the use of random probability measures to make inference on the distribution of the data. The randomness, however, is of different nature. Frequentist methods build estimators whose randomness is only due to the randomness of the data. Standard examples include, e.g., empirical processes, histograms, bootstrapping, kernel density estimators, splines, and wavelets (see, e.g, \citet{Wasserman2006} and references therein). In a Bayesian setting the distribution of the data  is treated as the unknown parameter, and thus modeled as a random probability $\tilde P$. The posterior is the distribution of $\tilde P$ \emph{conditionally} on the observed data, thus in principle its randomness does not stem from the data but from the uncertainty about the parameter. The most popular random probability in Bayesian nonparametrics is the Dirichlet process \cite{Ferguson1973}, though many useful generalizations are available in the literature, including species sampling processes, normalized completely random measures, stick-breaking processes, and kernel mixtures (see, e.g, \citet{Ghosal2017} and references therein). 

Several interesting theoretical and applied findings in nonparametric inference may be framed as approximations of some distribution of interest. These include approximation to a ground truth, as in the study of consistency, to the distribution of another population, as in two-sample tests, to another posterior, as in the merging of opinions, or to the exact but unattainable posterior, as in (Markov chain) Monte-Carlo or other finite-dimensional approximations. In most of these problems to quantify the quality of the approximation one needs a distance between laws of random probabilities, a.k.a., probability distributions on spaces of probability distributions. Such a distance should ideally satisfy the following properties:

\begin{enumerate}
\item {\em Metrization of weak convergence.} Many of the most common random probabilities, such as empirical processes and Dirichlet processes, are almost surely discrete and are typically mutually singular. Weak convergence provides meaningful comparisons between distributions with different supports. 

\item {\em Numerical estimation through samples.} Frequently in applications the exact distribution of the random probability is not known and thus the distance must be approximated numerically through samples, with both low sample complexity (the approximation error as the sample size increases) and low computational complexity (the speed of the algorithm as the sample size increases).

\item {\em Analytical upper and lower bounds.} To prove asymptotic statements and to perform theoretical comparisons with numerical estimations, one needs analytical expressions that capture the most meaningful behavior of the distance.
\end{enumerate}

A natural candidate is the Wasserstein distance on $\mathcal{P}(\mathcal{P}(\X))$, the space of laws of random probabilities on the sample space $\X$. Indeed, the Wasserstein distance on a Polish space $\X$ makes the space of probabilities $\mathcal{P}(\X)$ a Polish space (cfr. Remark 7.1.7 in \citet{ambrosio2005gradient}). 
Thus, one can build a Wasserstein distance on $\mathcal{P}(\mathcal{P}(\X))$ using the Wasserstein distance on $\mathcal{P}(\X)$ as ground metric. We refer to the induced metric as \emph{Wasserstein over Wasserstein distance}, also known under the name \emph{Hierarchical Optimal Transport} distance in the optimal transport community. This elegant and intuitive metric has been recently independently defined and used in many different contexts, including the analysis of convergence of the mixing measure in Bayesian hierarchical mixture models \cite{Nguyen2016}, to measure the similarity measure between documents in topic models \cite{Yurochkin2019, bing2022annals, bing2022sketched}, for joint clustering of observations and their distributions \cite{Ho2017}, as training loss in generative adversarial networks on images \cite{Dukler2019}, and to measure the discrepancy between datasets for classification tasks \cite{AlvarezMelis2020}.

In this work we show that though it is indeed possible to estimate the Wasserstein over Wasserstein distance numerically through samples, the infinite dimensionality of $\mathcal{P}(\X)$ drastically deteriorates its sample complexity, which can be slower than any polynomial rate of the sample size (Theorem~\ref{th:bad_rates_wow}). A classical workaround to the curse of dimensionality in optimal transport is to use entropic regularization \cite{mena2019statistical} but it is unclear if it would be enough here, as the base space $\mathcal{P}(\X)$ is infinite-dimensional. Since in practice a very large sample size may not be available or, more often, it may lead to an insurmountable computational burden, we propose a new distance between the laws of random probabilities that can achieve a parametric rate of convergence (Definition \ref{def:hipm} and Theorem \ref{th:hier_estimator}). One of the main ingredients we use are integral probabilities metrics \cite{Mueller1997}, and for this reason we call our distance \emph{hierarchical IPM}. The hierarchical IPM is dominated by the Wasserstein over Wasserstein distance, and thus retains all analytical upper bounds of the former. We also provide a general strategy to compute lower bounds and use them to recover closed-form expressions of the Wasserstein over Wasserstein distance when applied to empirical measures, to Dirichlet processes, and to the more general species sampling processes (Theorem \ref{th:species_sampling}). Moreover, in Theorem \ref{th:weak_convergence} we show that our distance metrizes the weak convergence on compact spaces. This completes our list of desirable properties for a distance on the laws of random probabilities. Note that the analysis we conduct here is valid for a fairly general class of IPMs, but the choice of the underlying function space affects the computational and statistical properties of our distance and may be chosen depending on its intended use.

The Dirichlet process is the cornerstone of Bayesian nonparametrics. It provides the law for a random probability that has been effectively used to build priors and derive posteriors for infinite-dimensional parameters (see \citet{Ghosal2017} and references therein). Its infinite dimensionality provides elegant and interpretable analytical properties that are less prone to computational algorithms. For this reason, a wealth of finite-dimensional approximations has been developed, each of which recovers the Dirichlet process as a limit. Whereas limiting behaviours are typically well-studied, approximation errors at the level of the random probability are rarely available. To our knowledge, these have been only carried out for almost-sure truncations of the random measures (see, e.g., \citet{MuliereTardella1998, IshwaranZarepour2000, IshwaranJames2001, Arbel2019}) but not for distributional approximations, where the error is usually studied in terms of the $L_1$ distance between the marginal densities of the data distribution induced by specific Bayesian models (see, e.g., \citet{IshwaranZarepour2000,IshwaranZarepour2002,Campbell2019,Lijoi2020, Nguyen2023}). As an application of our findings, we use both the Wasserstein over Wasserstein distance and our HIPM to investigate the quality of some of the most common finite-dimensional approximations, namely the Dirichlet multinomial process, the truncated stick-breaking, and the hierarchical empirical measure.
We provide numerical estimation of the approximation error and are able to compute non-asymptotic analytical error bounds that provide new rules-of-thumb for deciding the best one to use in practice.  

In summary, the main contributions of our work are i) drawing attention on the need for a distance between the laws of random probabilities; ii) showing that the sample complexity of the Wasserstein over Wasserstein distance $\W_\W$ can be slower than any polynomial rate in the sample size (Section \ref{sec:statistics}); iii) proposing a new distance, the hierarchical IPM, which is strongly related to $\W_\W$ but with better sample complexity (Sections \ref{sec:background} and \ref{sec:statistics}) ; iv) using our distance to provide closed-form expressions for $\W_\W$ on empirical measures and on species sampling processes (Section \ref{sec:topology}); iv) proposing a gradient ascent algorithm to compute our distance (Section \ref{sec:numerics}); v) applying our distance to find the best finite-dimensional approximation of the Dirichlet process in Bayesian nonparametrics (Section \ref{sec:application}).

\section{Background and main definition}
\label{sec:background}

Let $(\mathbb{X},d_\X)$ be a Polish (metric, complete and separable) space with a bounded diameter ${\rm diam}(\X) = \sup_{x,y \in \X} d_\X(x,y)$, e.g. a compact subset of $\R^d$ endowed with the Euclidean distance. We focus on probability distributions on the space $\mathcal{P}(\X)$ of probabilities on $\mathbb{X}$, denoted by $\mathcal{P}(\mathcal{P}(\mathbb{X}))$. Elements of $\mathcal{P}(\mathbb{X})$ will usually be denoted as $P$, whereas elements of  $\mathcal{P}(\mathcal{P}(\mathbb{X}))$ will usually be denoted by $\mathbb{Q}$. Random elements on these spaces are distinguished by a $\sim$, e.g. $\tilde P$ denotes a random probability. For any measurable function $f$ on $\X$, $P(f) = \int_\X f \, \ddr P$. For a $\X$-valued random variable $X$ we denote its law by $\mathcal{L}(X) \in \mathcal{P}(\X)$. We write $\eqd$ for equality in distribution of random variables, that is, $X \eqd Y$ if $\mathcal{L}(X) = \mathcal{L}(Y)$. Eventually, we denote by $\E$ the expectation of random variables, and we use $\E_X$ to emphasize the source of randomness.  

We define distances between the laws of random probabilities by using two baseline ingredients: integral probability metrics \cite{Zolotarev1984, Mueller1997} and the Wasserstein distance on a generic bounded Polish space $(\mathbb{Y},d_\Y)$, which in this work will be either $\X$ or the space of probabilities $\mathcal{P}(\X)$ with a suitable metric.  \\

\begin{definition}
Let $\mathcal{F}$ be a class of $\R$-valued bounded measurable functions on a Polish space $(\mathbb{Y},d_\Y)$. The integral probability metric (IPM) between $P_1,P_2 \in \mathcal{P}(\mathbb{Y})$ is 
\[
\mathcal{I}_{\mathcal{F}}(P_1,P_2) = \sup_{f \in \mathcal{F}} | P_1(f) - P_2(f)|.
\]
\end{definition}

A finite IPM is a distance whenever $\mathcal{F}$ separates probabilities on $\mathbb{Y}$, i.e. if $P_1(f) = P_2(f)$ for every $f \in \mathcal{F}$ implies that $P_1 = P_2$. IPMs encompass many well-established distances between probability measures, including the total variation distance, the Maximum Mean Discrepancy (MMD) for a characteristic kernel and the Wasserstein distance of order 1, which is recovered when $\mathcal{F} = {\rm Lip}_1(\mathbb{Y}, \R)$ is the class of 1-Lipschitz functions on $(\mathbb{Y},d_\Y)$ \citep{sriperumbudur2012}. By duality (see Remark 6.5 in \citep{Villani2008}) the Wasserstein distance can also be expressed as an infimum over all couplings as follows. We recall that a coupling between two probabilities $P_1, P_2$ on a Polish space $(\mathbb{Y},d_\Y)$ is the law of any random vector $(Y_1,Y_2)$ on the product space $\mathbb{Y} \times \mathbb{Y}$ such that $Y_1 \sim P_1$ and $Y_2 \sim P_2$. We denote by $\Gamma(P_1,P_2)$ the set of couplings between $P_1$ and $P_2$. \\

\begin{definition}
\label{def:wass_classical}
The Wasserstein distance of order $p$ between $P_1,P_2 \in \mathcal{P}(\mathbb{Y})$ is 
\[
\mathcal{W}_p(P_1, P_2)^p = \inf_{\gamma \in \Gamma(P^1, P^2)} \mathbb{E}_{(Y_1,Y_2) \sim \gamma}(d_\Y(Y_1,Y_2)^p).
\]
\end{definition}

We focus on $p=1$ because of its link with IPMs and for simplicity we denote the Wasserstein of order 1 on $\X$ as $\mathcal{W}$. When referring to the Wasserstein distance between the laws of two random variables $Y_1,Y_2$ we will sometimes omit the law $\mathcal{L}$ and write $\W(Y_1,Y_2) = \W(\mathcal{L}(Y_1), \mathcal{L}(Y_2))$.

When the space of probabilities $\mathcal{P}(\X)$ endowed with an IPM is a metric space, the Wasserstein distance over this IPM naturally defines a distance between the laws of random probabilities. For an IPM $\mathcal{I}_\mathcal{F}$ we denote by $\rm{Lip}_1(\mathcal{I}_\mathcal{F}) = \mathrm{Lip}_1((\mathcal{P}(\X), \mathcal{I}_\mathcal{F}),\R)$ the class of $\mathcal{I}_\mathcal{F}$-Lipschitz functions $h: \mathcal{P}(\mathbb{X}) \to \mathbb{R}$ s.t. $|h(P_1) - h(P_2)| \le \mathcal{I}_\mathcal{F}(P_1,P_2)$ for every $P_1, P_2 \in  \mathcal{P}(\mathbb{X})$. \\

\begin{definition}
\label{def:wow}
Let $\mathcal{I}_{\mathcal{F}}$ be an IPM on $\X$. Then for any $\mathbb{Q}_1, \mathbb{Q}_2 \in \mathcal{P}(\mathcal{P}(\X))$ we define
\begin{align*}
\W_{\mathcal{F}}(\mathbb{Q}_1, \mathbb{Q}_2) 
&= \inf_{ \gamma \in \Gamma(\mathbb{Q}_1, \mathbb{Q}_2)} \E_{(\tilde P_1, \tilde P_2) \sim \gamma} ( \mathcal{I}_{\mathcal{F}}(\tilde P_1, \tilde P_2)) \\ \nonumber
&= \sup_{h \in \rm{Lip}_1(\mathcal{I}_\mathcal{F})} | \mathbb Q_1(h) - \mathbb Q_2(h) |.
\end{align*}
We call $\W_{\mathcal{F}}$ the Wasserstein over IPM distance. \\
\end{definition}

For $\mathcal{F} = {\rm Lip}_1(\mathbb{X}, \R)$, $\rm{Lip}_1(\mathcal{I}_\mathcal{F}) = \rm{Lip}_1(\W)$ is the class of Wasserstein-Lipschitz functions  
on $\mathcal{P}(\mathbb{X})$. We thus recover the Wasserstein over Wasserstein distance, which we denote as $\W_{\W}$. As discussed in the introduction and investigated in Section~\ref{sec:statistics}, the statistical properties of $\W_\W$ are not satisfactory and for this reason we propose a new, yet related, distance. 
The main idea is that for a random probability $\tilde P$, the laws of the integrals $\tilde P(f)$, where $f$ is a continuous bounded function, are probabilities on $\R$ and they are enough to characterize the law of $\tilde P$ (cfr. Theorem 4.11 in \citet{Kallenberg2017}). \\

\begin{definition}
\label{def:hipm}
Let $\mathcal{I}_{\mathcal{F}}$ be an IPM on $\X$. Then for any $\mathbb{Q}_1, \mathbb{Q}_2 \in \mathcal{P}(\mathcal{P}(\X))$ we define 
\begin{equation*}
d_{\mathcal{F}}(\mathbb{Q}_1, \mathbb{Q}_2) = \sup_{f \in \mathcal{F}} \W ( \mathcal{L}(\tilde P_1(f)) , \mathcal{L}(\tilde P_2(f)) ),
\end{equation*}
where $ \tilde P_i \sim \mathbb{Q}_i$ for $i=1,2$. We call $d_{\mathcal{F}}$ the hierarchical integral probability metric (HIPM). \\
\end{definition}

Conditions guaranteeing that $d_{\mathcal{F}}$ is a distance, together with a characterization of the topology it generates, can be found below in Theorem~\ref{th:weak_convergence}.
In the following we will devote much attention to $\mathcal{F} = {\rm Lip}_1(\mathbb{X}, \R)$. The corresponding metric $d_{\mathcal{F}}$ will be compactly denoted $\dlip$ and we refer to it as the Lipschitz HIPM. 

Note that an analogy can be made with the max-sliced Wasserstein distance \citep{deshpande2019max} as $d_{\mathcal{F}}$ is 
the maximal one-dimensional Wasserstein distance between the ``projections'' of $\mathbb{Q}_1$, $\mathbb{Q}_2$ onto $\mathcal{P}(\R)$ via the maps $\mathbb{Q} \mapsto \mathcal{L}(\tilde{P}(f))$, with $\tilde{P} \sim \mathbb{Q}$. Another point of view is to see $d_\mathcal{F}$ as an IPM on $\mathcal{P}(\X)$: by expressing the Wasserstein distance on $\R$ as an IPM we observe that
\[
d_\mathcal{F}(\mathbb Q_1, \mathbb Q_2) = \sup_{h \in \mathfrak{F}} | \mathbb Q_1(h) - \mathbb Q_2(h) |,
\]
where here the class of functions $\mathfrak{F}$ from $\mathcal{P}(\X)$ to $\R$ is
\[
\mathfrak{F} = \{ P \mapsto g(P(f)) \ : \ f \in \mathcal{F} \text{ and } g \in \mathrm{Lip}_1(\R) \}.
\]
Thus both $\W_\mathcal{F}$ and $d_\mathcal{F}$ are special cases of IPMs on $\mathcal{P}(\X)$, whose corresponding classes of functions, $\rm{Lip}_1(\mathcal{I}_\mathcal{F})$ and $\mathfrak{F}$ respectively, satisfy $\mathfrak{F} \subset \rm{Lip}_1(\mathcal{I}_\mathcal{F})$. Intuitively, the class $\mathfrak{F}$ is much smaller than $\rm{Lip}_1(\mathcal{I}_\mathcal{F})$, and this is explains the better sample complexity of $\dlip$. This will be precisely quantified in terms of Rademacher complexity in Section~\ref{sec:statistics}. On the other hand, under reasonable assumptions both distances metrize the same topology, as discussed in Section~\ref{sec:topology}.

We conclude this section by recalling the paradigmatic example of a random probability: the Dirichlet process \cite{Ferguson1973}, which will be useful both in the study of the sample complexity of the Wasserstein over Wasserstein distance in Section~\ref{sec:statistics}, and in Section~\ref{sec:application}. We introduce it through its stick-breaking representation \cite{Sethuraman1994}. \\

\begin{definition}
\label{def:dp}
A random probability $\tilde P$ has a Dirichlet process distribution with concentration $\alpha>0$ and base probability $P_0 \in \mathcal{P}(\X)$, written $\tilde P \sim \textup{DP}(\alpha, P_0)$, if
\[
\tilde P \eqd \sum_{i=1}^{+\infty} J_i \delta_{X_i},
\]
where $X_i \simiid P_0$ are independent of $V_i \simiid \textup{Beta}(1,\alpha)$, and $J_i = V_i \prod_{j=1}^{i-1}(1- V_{j})$. The $\{J_i\}_i$ are termed the stick-breaking weights. 
\end{definition}

\section{Topological and metric properties}
\label{sec:topology}

In this section we prove that our distance metrizes the desired topology and explain its relation to the more-established Wasserstein over Wasserstein distance. We write $(C_b(\X), \mathrm{L}_\infty)$ for the space of continuous and bounded functions endowed with the supremum norm. We recall that a sequence $(P_n)_{n}$ of probability laws on a Polish space is said to converge \emph{weakly} to $P$ if for every $f \in C_b(\X)$, $P_n(f)$ converges pointwise to $P(f)$. This turns $\mathcal{P}(\X)$ into a Polish space, see Remark 7.1.7 in \citet{ambrosio2005gradient}. Thus we can endow $\mathcal{P}(\mathcal{P}(\X))$ with the weak convergence when $\mathcal{P}(\X)$ itself is endowed with weak convergence. Remarkably, one can prove (cfr. \citet{Prohorov1961, vonWaldenfels1968, Harris1971} and Theorem 4.11 in \citet{Kallenberg2017}) that a sequence $(\mathbb{Q}_n)_{n}$ of laws of random probabilities converges weakly to $\mathbb Q$ if and only if for every $f \in C_b(\X)$, the sequence of real-valued random variables $\tilde P_n(f)$, for $\tilde P_n \sim \mathbb{Q}_n$, converges weakly to $\tilde P(f)$, for $\tilde P \sim \mathbb Q$. \\

\begin{theorem}
\label{th:weak_convergence}
If $\mathcal{I}_\mathcal{F}$ is bounded and metrizes the weak convergence on $\mathcal{P}(\X)$ then $\W_\mathcal{F}$ metrizes the weak convergence on $\mathcal{P}(\mathcal{P}(\X))$. If in addition $\{ af +b \ : \ a,b \in \R, \, f \in \mathcal{F} \}$ is dense in  $(C_b(\X), \mathrm{L}_\infty)$ then $d_\mathcal{F}$ also metrizes the weak convergence on $\mathcal{P}(\mathcal{P}(\X))$.  \\
\end{theorem}

In the compact case the requirement of density can be checked with the following lemma. \\

\begin{lemma}
\label{th:density}
Assuming $\X$ compact, let $\mathcal{F}$ a class of functions which is absolutely convex (that is, convex and such that $a f \in \mathcal{F}$ for any $|a| \leq 1$ and $f \in \mathcal{F}$) and such that $\mathcal{I}_\mathcal{F}$ is a distance. Then $\{ af +b \ : \ a,b \in \R, \, f \in \mathcal{F} \}$ is dense in  $(C_b(\X), \mathrm{L}_\infty)$. \\ 
\end{lemma}

Note that we can always inflate $\mathcal{F}$ so that it becomes absolutely convex but the IPM $\mathcal{I}_\mathcal{F}$ stays the same: see Theorem 3.3 in \citet{Mueller1997}.
As examples, $\dlip$ metrizes weak convergence if $\X$ is compact and so does $d_\mathcal{F}$ if $\mathcal{I}_\mathcal{F}$ is a MMD distance for a $c$-universal kernel (by definition of universality, see \citet{sriperumbudur2011universality}).

We move on to the study of analytical upper and lower bounds of these distances. For a random probability $\tilde P$ we denote by $\mathbb{E}(\tilde P)$ the deterministic mean measure which satisfies $\mathbb{E}(\tilde P)(A) = \mathbb{E} (\tilde P(A))$ for any measurable set $A$. \\

\begin{proposition}
\label{th:bounds}
Let $\W_{\mathcal{F}}$ denote the Wasserstein over $\mathcal{I}_{\mathcal{F}}$ distance. Then if $\tilde P_i \sim \mathbb{Q}_i$ for $i=1,2$,
\[
\mathcal{I}_{\mathcal{F}}(\mathbb{E}(\tilde P_1),\mathbb{E}(\tilde P_2))  \le  d_{\mathcal{F}}(\mathbb{Q}_1, \mathbb{Q}_2) \le  \W_{\mathcal{F}}(\mathbb{Q}_1, \mathbb{Q}_2).
\] \\
\end{proposition}

It is easy to see that the inequalities are tight. Indeed, when $\mathbb{Q}_1 = \delta_{P_1}$, and $\mathbb{Q}_2 = \delta_{P_2}$ for some (deterministic)  probability measures $P_1,P_2$,
\[
\mathcal{I}_{\mathcal{F}}(P_1, P_2) = d_{\mathcal{F}}(\mathbb{Q}_1, \mathbb{Q}_2) =  \W_{\mathcal{F}}(\mathbb{Q}_1, \mathbb{Q}_2).
\] 
In particular, for $\mathcal{F} = \rm{Lip}_1(\X,\R)$ our distance $\dlip$ is always a lower bound to the Wasserstein over Wasserstein distance $\W_\W$. The next result provides a general condition for equality to hold. We recall that $T:\X \to \X$ is an optimal transport map if its graph is a $c$-cyclically monotone subset of $\X \times \X$, in the sense of Definition 5.1 in \citet{Villani2008}. Moreover, if $P \in \mathcal{P}(\X)$, $T_\#P(A) = P(T^{-1}(A))$ is the push-forward measure of $P$ by the map $T$. In the following we write $\W_\W( \tilde P_1, \tilde P_2)$ in place of $\W_\W( \mathcal{L}(\tilde P_1), \mathcal{L}(\tilde P_2))$. \\

\begin{lemma}
\label{th:equality}
Let $T : \X \mapsto \X$ be an optimal transport map. If $T_\# \tilde{P}_1 \eqd \tilde{P}_2$,
\[
\W_{\W}( \tilde{P}_1, \tilde{P}_2 )  = \dlip( \tilde{P}_1, \tilde{P}_2 ) = \W( \mathbb{E}(\tilde P_1),\mathbb{E}(\tilde P_2) ).
\]
\end{lemma}

We now state a corollary of Lemma~\ref{th:equality} that shows its far reach. Let $\tilde P_1$ and $\tilde P_2$ be two discrete random probabilities
\[
\tilde P_1 = \sum_{j \ge 1} J_j^{(1)} \delta_{X_j^{(1)}}, \qquad \tilde P_2 = \sum_{j \ge 1} J_j^{(2)} \delta_{X_j^{(2)}},
\]
such that (a) the jumps are independent of the atoms, i.e., $(J_j^{(i)})_{ j \ge 1} \perp (X_j^{(i)})_{j \ge 1;}$ for $i=1,2$; (b) the atoms are i.i.d. from base distributions $P_1, P_2$ respectively, i.e. $X_j^{(i)} \simiid P_i$ for $i=1,2$; (c) the jump distribution of $\tilde P_1$ and $\tilde P_2$ are the same, i.e. $(J_j^{(1)})_{ j \ge 1} \eqd (J_j^{(2)})_{ j \ge 1.}$. \\

\begin{theorem}
\label{th:species_sampling}
Let $\tilde P_1$ and $\tilde P_2$ be s.t. (a), (b), (c) hold. Then,
\[
\W_{\W}( \tilde{P}_1, \tilde{P}_2 )  = \dlip( \tilde{P}_1, \tilde{P}_2 ) = \W(P_1,P_2).
\]
\end{theorem}

Theorem~\ref{th:species_sampling} holds both for probabilities with a finite and an infinite number of atoms. In particular, the assumptions (a) and (b) are the ones that define the notable class of \emph{species sampling processes} \cite{Pitman1996}. This is a very general class of random probabilities which encompasses, e.g., empirical measures, the Dirichlet process, the Pitman-Yor process \cite{PitmanYor1997}, homogeneous normalized completely random measures \cite{Regazzini2003}: Theorem~\ref{th:species_sampling} applies to all these examples. The contribution is two-fold: on the one hand it shows that in many cases our distance coincides with $\W_\W$, making it a good alternative, on the other it provides the exact expression of $\W_\W$ in many interesting scenarios. We point out two useful corollaries. \\

\begin{corollary}
\label{th:empirical}
Let $X_1,\dots, X_n \simiid P_1$ and $Y_1,\dots, Y_n \simiid P_2$. Then,
\[
\W_{\W}\bigg( \mathcal{L} \left( \frac{1}{n} \sum_{i=1}^n \delta_{X_i} \right), \mathcal{L} \left( \frac{1}{n} \sum_{i=1}^n \delta_{Y_i} \right) \bigg) =  \W(P_1,P_2).
\]
\end{corollary}

Corollary~\ref{th:empirical} highlights the crucial role of the dependence across samples $\{X_i\}_i$ and $\{Y_i\}_i$ in estimating $\W(P_1,P_2)$. Whereas it is known that with independent samples the bias goes slowly to zero when $n$ increases and suffers the curse of dimensionality, considering perfectly coupled atoms removes the bias. We thus expect partial forms of dependence across samples to also reduce the bias. 

The next result focuses on the Dirichlet process (Definition~\ref{def:dp}).\\

\begin{corollary}
\label{th:exact_dp}
For $P_1, P_2 \in \mathcal{P}(\X)$  and $\alpha > 0$,
\[
\W_{\W}(\rm{DP}(\alpha, P_1), \rm{DP}(\alpha, P_2)) = \W(P_1,P_2).
\]
\end{corollary}
This recovers the prominent identity in \citet[Lemma 3.1]{Nguyen2016}. The author addresses it as a `remarkable identity of the Dirichlet process' - our result extends it to all species sampling models.

Another useful feature of $\mathcal{F} = \rm{Lip}_1(\R,\R)$ is that the identity function $f(x) = x$ belongs to $\mathcal{F}$. Then a natural lower bound to $\dlip$, and thus $\W_\W$, is the standard Wasserstein distance between the random means: for $\tilde P_i \sim \mathbb{Q}_i$,
\begin{equation}
\label{eq:lower}
\W \bigg(  \int_\R x \, \ddr \tilde P_1(x) ,  \int_\R x \, \ddr \tilde P_2(x) \bigg) \le \dlip ( \mathbb{Q}_1,  \mathbb{Q}_2).
\end{equation}

\section{Sample complexity}
\label{sec:statistics}

Let $\tilde P_1,\dots,\tilde P_n \simiid \mathbb{Q}$ and consider the empirical estimator 
\[
\mathbb{\tilde Q}_{(n)} = \frac{1}{n} \sum_{i=1}^n \delta_{\tilde P_i}.
\]
We have introduced a third level of randomness: $\mathbb{\tilde Q}_{(n)}$ is a random variable taking values in $\mathcal{P}(\mathcal{P}(\X))$, with the randomness coming from $\tilde P_1,\dots,\tilde P_n$.  
The extension of Glivenko-Cantelli theorem to general Polish spaces \cite{Varadarajan1958} ensures that $\mathbb{\tilde Q}_{(n)}$ converges weakly to $\mathbb{Q}$ almost surely (a.s.). By Theorem~\ref{th:weak_convergence}, both $\W_{\W}(\mathbb{\tilde Q}_{(n)}, \mathbb{Q})$ and $\dlip(\mathbb{Q}_{(n)}, \mathbb{Q})$ converge to zero a.s. as the sample size $n$ diverges.   
Yet, we prove that there is a crucial difference in their convergence rate: whereas the convergence of the former is slower than polynomial, the latter can achieve parametric rate of convergence. \\

\begin{theorem}
\label{th:bad_rates_wow}
Let $\X$ a bounded subset of $\R^d$. Then there exist constants $C$, $N$ depending on $d$ and $\rm{diam}(\X)$ such that, for $n \geq N$, for any $\mathbb{Q} \in \mathcal{P}(\mathcal{P}(\X))$
\[
\E(\W_\W(\mathbb{\tilde Q}_{(n)}, \mathbb{Q})) \le C \frac{ \log(\log(n))}{\log(n)}.
\]
Moreover let $\mathbb{Q} = {\rm DP}(\alpha, P_0)$ be the law of a Dirichlet process for $\alpha > 0$ and $P_0 \in \mathcal{P}(\X)$ whose support has a non-empty interior. Then for every $\gamma > 0$ there exists $c_\gamma > 0$ such that, at least for $n$ large enough,
\[
\E(\W_\W(\mathbb{\tilde Q}_{(n)}, \mathbb{Q})) \ge c_\gamma n^{-\gamma}.
\]
\end{theorem}

On the other hand we show that $d_{\mathcal{F}}$ has parametric convergence rate when $\X \subseteq \R$ under some uniformly bounded condition on $\mathcal{F}$ with respect to the supremum norm. We give a slightly weaker condition in terms of $\mathcal{F}^* = \{ f^* = f - f(x_0) \text{ s.t. } f \in \mathcal{F} \}$ for a fixed $x_0 \in \X$. This will allow one to treat the case $\mathcal{F} = \rm{Lip}_1(\X,\R)$, as $\rm{Lip}_1^*(\X,\R)$ is uniformly bounded whereas $\rm{Lip}_1(\X,\R)$ contains all constant functions. \\

\begin{lemma}
\label{th:entropy}
Let $\mathcal{F}^* = \{ f^* = f - f(x_0) \text{ s.t. } f \in \mathcal{F} \}$ be uniformly bounded by $K$ in the supremum norm, for a fixed $x_0 \in \X$. Then for every $\mathbb{Q} \in \mathcal{P}(\mathcal{P}(\X))$, 
\begin{equation*}
\E(d_{\mathcal{F}}(\mathbb{\tilde Q}_{(n)}, \mathbb Q)) \le \frac{320 \log(2) K}{\sqrt{n}} + \inf_{\epsilon>0} \bigg \{4 \epsilon + \frac{64}{\sqrt{n}} \int_{\epsilon/4}^{K} \sqrt{ \log N \left( \frac{\delta}{2} ; \mathcal{F}^*, \rm{L}_{\infty} \right)} \, \ddr \delta \bigg\},    
\end{equation*}
where $N(\cdot; \mathcal{F}^*, \rm{L}_{\infty})$ is the covering number of $\mathcal{F}^*$ with respect to the supremum norm. \\
\end{lemma}

Lemma~\ref{th:entropy} reduces the convergence rate of the HIPM to the convergence rate of the corresponding IPM (the infimum in the right hand side) as the first term coming from the ``hierarchical'' part of the distance goes to zero at the parametric rate $1/\sqrt{n}$. The infimum is already well-studied for a variety of classes, in particular for $\mathcal{F} = \rm{Lip}_1(\X, \R)$ when $\X$ is a subset of $\R^d$. \\

\begin{theorem}
\label{th:lipschitz_rates} 
Let $\mathcal{F} = \rm{Lip}_1(\X,\R)$ with $\X $ a bounded subset of $\R^d$. Then there exists a constant $C_d>0$ depending on $d$ and $\rm{diam}(\X)$ but not on $n$, such that, at least for $n$ large enough, 
\[
\E(\dlip(\mathbb{\tilde Q}_{(n)}, \mathbb{Q})) \le 
\begin{cases}
C_1 n^{-1/2} \qquad &\text{if } d = 1,\\
C_2 n^{-1/2} \log(n) \qquad &\text{if } d = 2,\\
C_d n^{-1/d} \qquad &\text{if }  d >2.
\end{cases} 
\]
\end{theorem}

In practice, to compute the HIPM distance between two empirical distributions we need that the realizations $\tilde P_1,\dots,\tilde P_n \simiid \mathbb{Q}$ are almost surely discrete and with a finite number of atoms. When this is not the case we may approximate $\tilde P_i$ through the empirical distribution of an exchangeable sequence whose de Finetti measure is $\tilde P_i$, namely $X_1^{(i)},\dots,X_m^{(i)} | \tilde P_i \simiid \tilde P_i$. By de Finetti's theorem 
$\tilde P_{i,(m)} = \frac{1}{m} \sum_{j=1}^m \delta_{X_j^{(i)}} \to \tilde P_i$
weakly almost surely. A hierarchical empirical estimator can then be defined as
\begin{equation}
\label{def:hier}
\mathbb{\tilde Q}_{(n,m)} = \frac{1}{n} \sum_{i=1}^n \delta_{\tilde P_{i,(m)}}.
\end{equation}
\begin{theorem}
\label{th:hier_estimator}
Let $\mathcal{F} = \rm{Lip}_1(\X,\R)$ with $\X $ a bounded subset of $\R^d$. Then there exist a constant $C_d>0$, not depending on $n$ or $m$, such that $\E(\dlip(\mathbb{\tilde Q}_{(n,m)}, \mathbb{Q}))$ is smaller or equal to
\[
\begin{cases}
C_1 (n^{-1/2} + m^{-1/2}) \qquad &\text{if } d = 1,\\
C_2 (n^{-1/2} \log(n) + m^{-1/2} \log(m)) \qquad &\text{if } d = 2,\\
C_d (n^{-1/d} + m^{-1/d}) \qquad &\text{if }  d >2.
\end{cases} 
\]
\end{theorem}

\section{Numerical estimation of the distances}
\label{sec:numerics}

In this section we focus on random measures that are fully discrete. That is, we consider $\mathbb{Q}_1, \mathbb{Q}_2 \in \mathcal{P}(\mathcal{P}(\X))$ such that
\[ 
\mathbb{Q}_k = \frac{1}{n} \sum_{i=1}^n \delta_{P^k_{i}},
\]
for $k=1,2$. Here $n$ is the number of distinct probabilities in the support of $\mathbb Q_1$, $\mathbb Q_2$. Each of the probabilities $P^k_i$ is also discrete and reads
\begin{equation}
\label{eq:discretization_lagrangian}
P^k_i = \frac{1}{m} \sum_{i=1}^m \omega_{i,j}^k \delta_{X^{k}_{i,j}} \qquad \text{with } \sum_{j=1}^{m} \omega^k_{i,j} = 1.
\end{equation} 

Thus $m$ is the number of atoms in each probability and we allow for non-uniform weights $\omega_{i,j}^k$ for each atom as it will be useful in Section~\ref{sec:application}. The measures $\mathbb{Q}_1$, $\mathbb{Q}_2$ can each be summarized as a $n \times m \times 2$ array containing the atoms $X^k_{i,j}$ and the weights $\omega^k_{i,j}$. Extensions to different values of $n$ for $\mathbb Q_1$ and $\mathbb Q_2$ and to a number $m$ of atoms which may depend on $i \in \{ 1, \ldots, n \}$ are possible, but we stick to this setting for the sake of simplicity. In this section we discuss how to output (an approximation of) $\W_\W(\mathbb{Q}_1,\mathbb{Q}_2)$ or $\dlip(\mathbb{Q}_1,\mathbb{Q}_2)$. 

\paragraph{Wasserstein over Wasserstein}

We follow Definition~\ref{def:wow}: we first need to compute the pairwise Wasserstein distance between $P^1_{i_1}$ and $P^2_{i_2}$ for any $i_1, i_2$, and then solve the $n \times n$ optimal transport problem whose cost matrix is given by the pairwise Wasserstein distances. This was already used e.g. in \citet{Yurochkin2019} and we refer to Appendix~\ref{app:algorithm} for the pseudocode. The only requirement is a package to solve the optimal transport problem. Assuming that the space $\X$ is a subset of $\R$, so that the inner optimal transport problems are just sorting problems, solving the $n \times n$ optimal transport exactly 
 requires $O(n^3 \log(n))$ operations \citep{orlin1997polynomial}. Thus this algorithm outputs the exact distance in $O(n^2 m \log(m) + n^3 \log(n))$ operations.

\paragraph{Lipschitz HIPM in dimension one}

We approximate our new distance $\dlip(\mathbb Q_1, \mathbb Q_2)$, as in Definition~\ref{def:hipm} with $\mathcal{F} = \mathrm{Lip}_1(\X,\R)$. We need to find a supremum over the space of Lipschitz functions $\mathrm{Lip}_1(\X,\R)$, and we will resort to a gradient ascent algorithm. We restrict $\X$ to be one-dimensional, so that $\X = [a,b]$. Higher dimensions would require a more careful parametrization of $\mathrm{Lip}_1(\X,\R)$, e.g. with neural networks, but we leave this investigation for future work.

To evaluate $P_i^k(f)$ for $k=1,2$ and $i=1, \ldots, n$, we only need the values of $f$ on $\{X^{k}_{i,j}\}$, that is on $2nm$ points. The requirement $f \in \mathrm{Lip}(\X,\R)$ corresponds to $2nm-1$ constraints of the form $|f(X) - f(X')| \leq |X-X'|$ for consecutive atoms $X,X'$. As we will run a gradient ascent with the function $f$ as unknown, this can become quickly prohibitive. For this reason, we will rather project the atoms of each measure on a fixed grid of stepsize $\Delta x$. Let $\ngr$ be the number of grid points and $Y_1, \ldots, Y_{\ngr}$ be the grid points over $[a,b]$. We consider $P^k_i$ given by 
\begin{equation}
\label{eq:discretization_Eulerian}
P^k_i = \sum_{q=1}^{\ngr} \omega^k_{i,q} \delta_{Y_q},  \qquad \text{with } \sum_{q=1}^{\ngr} \omega^k_{i,q} = 1,     
\end{equation}
which is a specific instance of~\eqref{eq:discretization_lagrangian} with common atoms. Each  random measure is now described by a $n \ngr$ array of weights $\omega^1_{i,q}$ and $\omega^2_{i,q}$. If $P^k_i$ is given by the form~\eqref{eq:discretization_lagrangian}, by projecting each atom $X^k_{i,j}$ onto the closest point on the grid $Y_1, \ldots, Y_{\ngr}$, it is possible to approximate it with a probability measure of the form~\eqref{eq:discretization_Eulerian}, up to an error $\Delta x$ in Wasserstein distance. As we would anyway have a statistical error of size $n^{-1/2} + m^{-1/2}$ if we are in the setting of Section~\ref{sec:statistics} (see Theorem~\ref{th:hier_estimator}), the rationale is that adding an error $\Delta x$ is reasonable if $\Delta x \asymp n^{-1/2} + m^{-1/2}$, that is, $\ngr \asymp \min(n^{1/2},m^{1/2})$.

We now expand Definition~\ref{def:hipm} in this setting. By Birkhoff theorem we rewrite the optimal transport problem as an infimum over permutations, and the optimal permutation can be found very efficiently by sorting; see Remark 2.2.27 in \citet{peyre2019computational}. We replace the function $f \in \mathrm{Lip}_1([a,b],\R)$ by a vector $\mathbf{f} \in \R^{\ngr}$ corresponding to $\mathbf{f}_q = f(Y_q)$ the evaluation of $f$ on grid points. The distance can be rewritten as 
\begin{align*}
\sup_{\mathbf{f} \in \R^{\ngr}} & \; \mathcal{G}(\mathbf{f})  \\
\text{such that} & \; |\mathbf{f}_{q+1} - \mathbf{f}_{q}| \leq \Delta x \; \forall q \in \{ 1, \ldots, {\ngr}-1 \}, 
\end{align*}
 where, by denoting $\mathcal{S}(n)$ the  permutations of $\{1,\ldots,n\}$,
\[
\mathcal{G}(\mathbf{f}) :=  \inf_{\sigma \in \mathcal{S}(n)} \frac{1}{n} \sum_{i=1}^n \left| \sum_{q=1}^{\ngr} (\omega^1_{i,q} - \omega^2_{\sigma(i),q}) \mathbf{f}_q  \right|.
\]

The constraints correspond to $\mathbf{f}$ being the restriction of a function in $\mathrm{Lip}_1(\X,\R)$ and they define a convex set. The function $\mathcal{G}$ we maximize is piece-wise linear. Being an infimum of convex functions, it is neither convex nor concave: we have left the realm of convex optimization. However, the gradient can be easily found, as explained in Appendix~\ref{app:algorithm}. Thus we use a gradient ascent algorithm to perform maximization of $\mathcal{G}$, but we have no guarantee of finding a global maximizer (in practice, we start with different random initializations). As the function is piece-wise linear there is no canonical stepsize for the gradient ascent, thus we implemented the following: starting from the gradient as an ascent direction, we first project it orthogonally on the set of ascent directions that preserve the convex constraint. Along this new direction, we perform a line search to move enough for the function to increase, and we stop if the increase is too small. We also do a linear change of variables and rather parametrize the function $\mathbf f$ by its derivative $\mathbf{g}$ (that is $\mathbf{f}_q = \sum_{q' < q} \Delta x \, \mathbf{g}_{q'}$ for any $q$) to simplify the Lipschitz constraint in a box constraint. The details and pseudocode can be found in Appendix~\ref{app:algorithm}.

\paragraph{Experiments}

\begin{figure}
\begin{center}

\begin{tabular}{cc}

\begin{tikzpicture}

\begin{axis}[
width=0.48\textwidth,
ymin = 0.35, ymax = 0.7,
xlabel={$n$},
ylabel={distance},
legend pos= south east]

\addplot[color=black, line width = 0.5pt, style = dashed] table [x=n, y=value]{graph_split_gt.txt};
\addlegendentry{{\small true value}}

\addplot[color=blue, line width = 0.5pt, mark = *] table [x=n, y=mean]{graph_split_nd.txt};
\addlegendentry{{\small $d_\mathrm{Lip}$}}

\addplot[color=red, line width = 0.5pt, mark = *] table [x=n, y=mean]{graph_split_wow.txt};
\addlegendentry{{\small $\W_\W$}}

\addplot[color=black!60!green, line width = 0.5pt, mark = *] table [x=n, y=mean]{graph_split_lb.txt};
\addlegendentry{{\small Bound~\eqref{eq:lower}}}

\addplot [name path=upper,draw=none] table[x=n,y expr=\thisrow{mean}+\thisrow{error}] {graph_split_nd.txt};
\addplot [name path=lower,draw=none] table[x=n,y expr=\thisrow{mean}-\thisrow{error}] {graph_split_nd.txt};
\addplot [fill=blue!10] fill between[of=upper and lower];

\addplot [name path=upper,draw=none] table[x=n,y expr=\thisrow{mean}+\thisrow{error}] {graph_split_wow.txt};
\addplot [name path=lower,draw=none] table[x=n,y expr=\thisrow{mean}-\thisrow{error}] {graph_split_wow.txt};
\addplot [fill=red!10] fill between[of=upper and lower];

\addplot [name path=upper,draw=none] table[x=n,y expr=\thisrow{mean}+\thisrow{error}] {graph_split_lb.txt};
\addplot [name path=lower,draw=none] table[x=n,y expr=\thisrow{mean}-\thisrow{error}] {graph_split_lb.txt};
\addplot [fill=black!10!green!10] fill between[of=upper and lower];

\end{axis}
\end{tikzpicture}

&


\begin{tikzpicture}
\begin{loglogaxis}[
width=0.48\textwidth,
xlabel={$n$},
ylabel={distance},
legend pos=south west]

\addplot[color=blue, line width = 0.5pt, mark = *] table [x=n, y=mean]{graph_same_nd.txt};
\addlegendentry{{\small$\dlip$}}

\addplot[color=red, line width = 0.5pt, mark = *] table [x=n, y=mean]{graph_same_wow.txt};
\addlegendentry{{\small $\W_\W$}}

\addplot[color=black!60!green, line width = 0.5pt, mark = *] table [x=n, y=mean]{graph_same_lb.txt};
\addlegendentry{{\small Bound~\eqref{eq:lower}}}

\addplot [name path=upper,draw=none] table[x=n,y expr=\thisrow{mean}+\thisrow{error}] {graph_same_lb.txt};
\addplot [name path=lower,draw=none] table[x=n,y expr=\thisrow{mean}-\thisrow{error}] {graph_same_lb.txt};
\addplot [fill=black!10!green!10] fill between[of=upper and lower];

\addplot [name path=upper,draw=none] table[x=n,y expr=\thisrow{mean}+\thisrow{error}] {graph_same_nd.txt};
\addplot [name path=lower,draw=none] table[x=n,y expr=\thisrow{mean}-\thisrow{error}] {graph_same_nd.txt};
\addplot [fill=blue!10] fill between[of=upper and lower];

\addplot [name path=upper,draw=none] table[x=n,y expr=\thisrow{mean}+\thisrow{error}] {graph_same_wow.txt};
\addplot [name path=lower,draw=none] table[x=n,y expr=\thisrow{mean}-\thisrow{error}] {graph_same_wow.txt};
\addplot [fill=red!10] fill between[of=upper and lower];

\end{loglogaxis}

\end{tikzpicture}  

\end{tabular}

\end{center}
\caption{Distances between independent realizations of $\tilde{\mathbb{Q}}_{(n,m)}$ and $\tilde{\mathbb{Q}}'_{(n,m)}$ in \eqref{def:hier}. We fix $m=5000$ and let $n$ evolve. Computations are repeated over $24$ realizations, with errors corresponding to one standard deviation. In the left plot $\mathbb{Q} = \rm{DP}(1,P_1)$ and $\mathbb{Q}'= \rm{DP}(1,P_2)$ with $P_1 = \mathrm{Unif}([-1/2,1/2])$, $P_2 = 1/2 \, \mathrm{Unif}([-1, -3/4]) + 1/2 \, \mathrm{Unif}([3/4, 1])$ so that $\W_\W(\mathbb{Q}, \mathbb{Q}') = \W(P_1,P_2) = 5/8$. In the right plot (in log-log scale) $\mathbb Q = \mathbb Q' = \text{DP}(1, \text{Unif}[0,1])$.}
\label{fig:distance_numerics}
\end{figure}
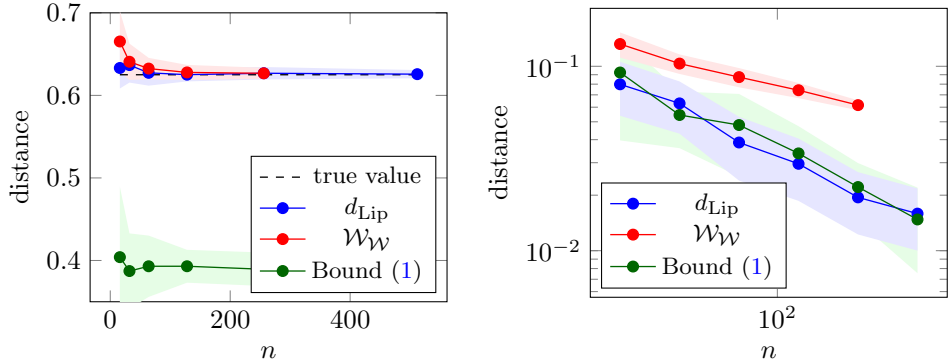

To illustrate our code and the results of Section~\ref{sec:statistics} we consider the setting where the discrete measures $\mathbb Q_1$ and $\mathbb Q_2$ are realizations of the empirical hierarchical estimators $\tilde{\mathbb Q}_{(n,m)}$, $\tilde{\mathbb Q}'_{(n,m)}$ in \eqref{def:hier} of two Dirichlet processes, $\mathbb Q$ and $\mathbb Q'$ respectively.

We first consider $\mathbb Q = \mathrm{DP}(P_1,\alpha)$ and $\mathbb Q' = \mathrm{DP}(P_2,\alpha)$ so that by Corollary~\ref{th:exact_dp} $\W_\W(\mathbb Q, \mathbb Q') = \dlip(\mathbb Q, \mathbb Q') = \W(P_1,P_2)$. We fix $m$ large and check that, as $n$ increases, the distances between realizations of $\tilde{\mathbb Q}_{(n,m)}$ and $\tilde{\mathbb Q}'_{(n,m)}$ indeed converge towards the true value. We also report the lower bound~\eqref{eq:lower} for completeness, and we observe that in this scenario it is informative but it is not tight. The results are displayed in the left plot of Figure~\ref{fig:distance_numerics}. 

To further emphasize the difference between $\W_\W$ and $\dlip$, we also consider the case where $\mathbb Q_1$ and $\mathbb Q_2$ are independent realizations of the same $\tilde{\mathbb Q}_{(n,m)}$, where $\mathbb Q$ is taken to be a Dirichlet process. We expect $\W_\W(\mathbb Q_1, \mathbb Q_2)$ and $\dlip(\mathbb Q_1, \mathbb Q_2)$ to go to zero. We fix $m$ large and study the convergence as $n \to + \infty$. We see that convergence is faster for $\dlip$ than for $\W_\W$ as predicted by our theory of Section~\ref{sec:statistics}. We also report the lower bound~\eqref{eq:lower}, and observe that it almost coincides with $\dlip$ in this case. The results are reported in the right plot of Figure~\ref{fig:distance_numerics}.

The algorithms were implemented in Julia. 
Computations were performed on the CPU of a standard 6-core laptop and we report execution time in the appendix. As computations were repeated on multiple realizations of $\tilde{\mathbb{Q}}_{n,m}$ and $\tilde{\mathbb{Q}}_{n,m}'$, we only keep cases where the distance for a single realization can be computed in less than one minute. This explains why some values of $n$ are missing for $\W_\W$ in Figure~\ref{fig:distance_numerics}. The code is available at the following address: \\
\begin{center}
\url{https://github.com/HugoLav/HierarchicalIPM}    
\end{center}

\section{Application: Approximation errors for the Dirichlet process}
\label{sec:application}

We consider three popular finite-dimensional approximations of the Dirichlet process (Definition~\ref{def:dp}). Each approximation is indexed by a finite number of atoms $N$ and, as $N$ goes to $+\infty$, it converges to a $\rm{DP}(\alpha, P_0)$.

\begin{enumerate}
\item Dirichlet multinomial process
\[
\tilde P_1|X_1,\dots,X_{N} \sim \rm{DP}\bigg(\alpha, \frac{1}{N} \sum_{i=1}^N \delta_{X_i}\bigg), \quad X_i \simiid P_0.
\]
It is equivalent to $\tilde P_1 = \sum_{i=1}^N J_i \delta_{X_i}$ where the jumps $(J_1, \ldots, J_N) \in \R^N$ follow a Dirichlet distribution of parameter $(\alpha/N, \ldots, \alpha/N)$ while the atoms are i.i.d. with law $P_0$.
\item Truncated stick-breaking process 
\[
\tilde P_2 = \sum_{i=1}^{N-1} J_i \delta_{X_i} + \bigg(1-\sum_{i=1}^{N-1} J_i\bigg) \delta_{X_N},
\]
where $X_i \simiid P_0$ and $J_1,\dots, J_{N-1}$ are the first $N-1$ stick-breaking weights (see Definition~\ref{def:dp}).
\item The dependent or hierarchical empirical measure
\[
\tilde P_3 = \frac{1}{N}\sum_{i=1}^N \delta_{Y_i},
\]
where $Y_1,\dots, Y_N | \tilde P \simiid \tilde P$ and $\tilde P\sim \rm{DP}(\alpha,P_0)$. We already considered it in Section~\ref{sec:statistics}.
\end{enumerate}

The empirical measure is less established in the Bayesian nonparametric literature as a finite-dimensional approximation of the Dirichlet process. Its atoms are exchangeable but not independent and they can be sampled through the Pólya urn or Chinese restaurant process scheme \cite{Blackwell1973}. In the following we use both analytical and empirical arguments to show that in many regimes it provides similar approximation errors to the widely used multinomial Dirichlet process, which has been independently defined by many different authors (see Section 4 of \citet{IshwaranZarepour2002}).
We state our results in terms of upper bounds of $\W_\W$ which, thanks to Proposition~\ref{th:bounds}, are also upper bounds of our distance $\dlip$. \\

\begin{proposition} 
\label{th:approx}
If $\tilde P \sim \rm{DP}(\alpha, P_0)$ with $P_0 \in \mathcal{P}(\X)$,
\begin{align*}
\W_\W( \tilde P, \tilde P_1) &\le \E\bigg(\W \bigg( P_0, \frac{1}{N}\sum_{i=1}^N \delta_{X_i} \bigg)\bigg), \\
\W_\W(\tilde P, \tilde P_2) &\le \E(d_\X(X_1, X_2))  \bigg( \frac{\alpha}{\alpha+1} \bigg)^N, \\
\W_\W(\tilde P, \tilde P_3) &\le \E_{\tilde P} \bigg( \W \bigg(\tilde P, \frac{1}{N}\sum_{i=1}^N \delta_{Y_i} \bigg) \bigg), 
\end{align*}
where $X_1, \ldots, X_N \simiid P_0$ and $Y_1,\dots, Y_N| \tilde P \simiid \tilde P$. 
In particular if $P_0 \in \mathcal{P}(\R)$ and $F_0$ denotes its c.d.f.,
\begin{align*}
\W_\W(\tilde P, \tilde P_1)  &\le \frac{1}{\sqrt{N}} \int_{-\infty}^{+\infty} \sqrt{F_0(x)(1- F_0(x))} \ddr x, \\
\W_\W(\tilde P, \tilde P_2) &\le 2 \bigg( \frac{\alpha}{\alpha+1} \bigg)^N \int_{-\infty}^{+\infty} F_0(x)(1- F_0(x)) \ddr x , \\
\W_\W(\tilde P, \tilde P_3) &\le \sqrt{\frac{\alpha}{N(\alpha+1)}} \int_{-\infty}^{+\infty} \sqrt{F_0(x)(1- F_0(x))} \ddr x.
\end{align*}
\end{proposition}

Proposition~\ref{th:approx} sheds light on a number of interesting properties. The Dirichlet multinomial process has polynomial convergence rate in the number of atoms $N$, whereas the stick-breaking has exponential convergence rate. However, this rate depends on $\alpha$ and for a diverging sequence of $\alpha$ it may fail to converge. The exponential convergence of the stick-breaking approximation can also be captured by moment summaries of the truncation error, as described  e.g. in \citet{IshwaranZarepour2000}, where it is also underlined that ``there appears to be no simple method for assessing the adequacy" of the Dirichlet multinomial approximation. Moreover, the potentially critical regime for the stick-breaking ($\alpha \to +\infty$) is rarely mentioned in the literature. This is not an uncommon situation 
since if we consider $n$ conditionally i.i.d. observations with de Finetti measure $\tilde P$, the posterior is a Dirichlet process with concentration parameter $\alpha + n$, which thus diverges $\asymp n$.
In particular, when approximating the posterior arising from a sequence of conditionally i.i.d. observations with a Dirichlet process prior one should always choose $N$ diverging faster than the number of observations $n$, which can be a problem with large sample sizes. On the contrary, the approximation of the Dirichlet multinomial and of the empirical measures do not show this high dependence on $\alpha$ but they are negatively affected by an increasing dimension of the space. Interestingly, the upper bounds for the Dirichlet multinomial and the hierarchical empirical measure are quite similar, and this similarity is also confirmed by the simulations.

\paragraph{Practical conclusions}

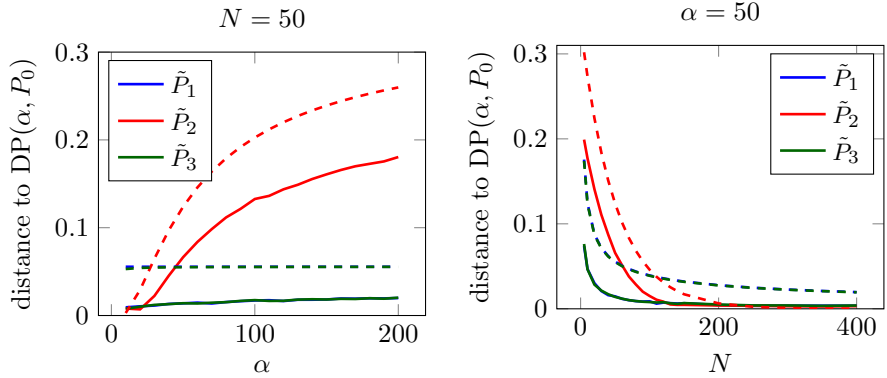
\begin{figure}
\begin{center}

\begin{tabular}{cc}

\begin{tikzpicture}
\begin{axis}[
width=0.45\textwidth,
ymin = 0.0, ymax = 0.3,
title = {$N = 50$},
xlabel={$\alpha$},
ylabel={distance to $\mathrm{DP}(\alpha,P_0)$},
legend pos=north west]

\addplot[color=blue, line width = 1pt] table [x=alpha, y=P1nd]{application_alpha.txt};
\addlegendentry{{\small $\tilde{P}_1$}}

\addplot[color=red, line width = 1pt] table [x=alpha, y=P2nd]{application_alpha.txt};
\addlegendentry{{\small $\tilde{P}_2$}}

\addplot[color=black!60!green, line width = 1pt] table [x=alpha, y=P3nd]{application_alpha.txt};
\addlegendentry{{\small $\tilde{P}_3$}}

\addplot[color=blue, line width = 1pt, style = dashed] table [x=alpha, y=P1ub]{application_alpha.txt};

\addplot[color=red, line width = 1pt, style = dashed] table [x=alpha, y=P2ub]{application_alpha.txt};

\addplot[color=black!60!green, line width = 1pt, style = dashed] table [x=alpha, y=P3ub]{application_alpha.txt};

\end{axis}
\end{tikzpicture}

&

\begin{tikzpicture}
\begin{axis}[
width=0.45\textwidth,
ymin = 0.0, ymax = 0.31,
title = {$\alpha = 50$},
xlabel={$N$},
ylabel={distance to $\mathrm{DP}(\alpha,P_0)$},
legend pos=north east]

\addplot[color=blue, line width = 1pt] table [x=N, y=P1nd]{application_N.txt};
\addlegendentry{{\small $\tilde{P}_1$}}

\addplot[color=red, line width = 1pt] table [x=N, y=P2nd]{application_N.txt};
\addlegendentry{{\small $\tilde{P}_2$}}

\addplot[color=black!60!green, line width = 1pt] table [x=N, y=P3nd]{application_N.txt};
\addlegendentry{{\small $\tilde{P}_3$}}

\addplot[color=blue, line width = 1pt, style = dashed] table [x=N, y=P1ub]{application_N.txt};

\addplot[color=red, line width = 1pt, style = dashed] table [x=N, y=P2ub]{application_N.txt};

\addplot[color=black!60!green, line width = 1pt, style = dashed] table [x=N, y=P3ub]{application_N.txt};

\end{axis}
\end{tikzpicture}
\end{tabular}
\end{center}
\caption{Distance from a $\rm{DP}(\alpha,P_0)$ for finite-dimensional approximations with $N$ atoms: Dirichlet multinomial ($\tilde P_1$), truncated stick breaking ($\tilde P_2$), and hierarchical empirical measure ($\tilde P_3$). The base measure $P_0$ is uniform over $[0,1]$. In the left plot $N=50$ is fixed, and $\alpha$ varies; in the right plot $\alpha=50$ is fixed, and $N$ varies. The solid line is the distance $\dlip$, while the dashed lines are the upper bounds of Proposition~\ref{th:approx}.}
\label{fig:finite_dimensional}
\end{figure}

Our study led to the following practical advice when using an approximation of the Dirichlet process. If one is approximating a Dirichlet processes with large concentration parameter,  as in the case of posteriors for a large dataset arising from a conditionally i.i.d. model directed by a Dirichlet process, the Dirichlet multinomial process and the hierarchical empirical measure will provide similar and more reliable approximations than the truncated stick-breaking, as depicted in the left plot of Figure~\ref{fig:finite_dimensional}. On the other hand if the observations live in a high-dimensional space or if we can afford taking $N$ consistently larger than the concentration parameter, the stick-breaking will tend to provide better approximations, as depicted in the right plot of Figure~\ref{fig:finite_dimensional}.

\section*{Limitations}

Our work presents some limitations from the computational and theoretical point of view. Our numerical method is restricted to a one-dimensional space and one should be careful in drawing delicate conclusions based on the numerical evaluation of the distance, as it relies on gradient ascent over a non-concave objective and thus it does not come with a guarantee of convergence. From a theoretical perspective, the sample complexity is derived for probabilities on a space with a bounded metric. Finally, we still lack a good understanding of the qualitative difference between the Wasserstein over Wasserstein distance and the HIPM.

\section*{Acknowledgments}

The authors thank Stefan Schrott and Giacomo Sodini for some helpful discussions as well as the anonymous referees for their suggestions which contrinute to improve the quality of the paper. The second author is also affiliated to the Bocconi Institute for Data Science and Analytics (BIDSA) and is partially supported by the MUR-Prin 2022-202244A7YL funded by the European Union - Next Generation EU.

\appendix 

\section{Proofs}

\subsection*{Proof of Theorem~\ref{th:weak_convergence}}

\emph{The distance $\W_\mathcal{F}$ metrizes weak convergence}. That $\W_\mathcal{F}$ metrizes weak convergence over $\mathcal{P}(\mathcal{P}(\X))$ is a simple consequence of Remark 7.1.7 in \citet{ambrosio2005gradient} and the fact that $(\mathcal{P}(\X), \mathcal{I}_\mathcal{F})$ has a bounded diameter by assumption.

\emph{The distance $d_\mathcal{F}$ metrizes weak convergence}. Let $\mathbb Q_n$ a sequence in $\mathcal{P}(\mathcal{P}(\X))$, and $\mathbb{Q} \in \mathcal{P}(\mathcal{P}(\X))$. We also write $\tilde P_n, \tilde P$ for random measures with $\tilde P_n \sim \mathbb Q_n$ and $\tilde P \sim \mathbb Q$. 

We first assume that $d_\mathcal{F}(\mathbb Q_n, \mathbb Q)$ converges to $0$ as $n \to + \infty$. In particular for any $g \in \mathcal{F}$, we know that $\W(\tilde P_n(g), \tilde P(g))$ converges to zero. This easily extends to any $g$ in $\mathcal{F}' = \{ af+b \ : \ a,b \in \R, \, f \in \mathcal{F} \}$. We extend to $C_b(\X)$ by our density assumption. If $f \in C_b(\X)$ is any continuous and bounded function and $\varepsilon > 0$, by assumption we can find $g \in \mathcal{F}'$ with $\| f- g \|_\infty \leq \varepsilon$. With the triangle inequality we easily obtain 
\begin{multline*}
\W(\tilde P_n(f), \tilde P(f)) \leq \W(\tilde P_n(f), \tilde P_n(g)) + \W(\tilde P_n(g), \tilde P(g)) + \W(\tilde P(g), \tilde P(f)) \\
 \leq  \W(\tilde P_n(g),\tilde P(g)) + 2 \varepsilon.
\end{multline*}
As $n \to + \infty$, the first term in the right hand side converges to $0$. As $\varepsilon$ can then be chosen arbitrary we have proved that $\W(\tilde P_n(f),\tilde P(f))$ converges to $0$ as $n \to + \infty$ for any continuous and bounded function $f$. Convergence in $\W$ implies converges in law, thus the random variable $\tilde P_n(f)$ converges in law to $\tilde P(f)$. This concludes the proof of weak convergence of $\mathbb Q_n$ to $\mathbb Q$. 

Conversely assume the weak convergence of $\mathbb Q_n$ to $\mathbb Q$. Using Proposition~\ref{th:bounds}, we know that $d_\mathcal{F}(\mathbb Q_n, \mathbb Q) \leq \W_\mathcal{F}(\mathbb Q_n, \mathbb Q)$. The first part of the present theorem yields that the right hand side goes to zero, and thus so does the left hand side.

\subsection*{Proof of Lemma~\ref{th:density}}

The proof strategy relies on the Hahn-Banach theorem and is similar to the one used to prove density results in the context of MMD distances, see e.g. Proposition 2 in \citet{sriperumbudur2011universality}.

\emph{Step 1. Characterization as linear span.} Using the absolute convexity of $\mathcal{F}$, we prove that $\mathcal{F}' = \{ af +b \ : \ a,b \in \R, \, f \in \mathcal{F} \}$ coincides with the linear span of $\mathcal{F} \cup \{ 1 \}$, where $1$ denotes the constant function equal to $1$. First we observe that $\mathcal{F}'$ is a vector space. Indeed, for any $a_1, a_2 \in \R$ and $f_1, f_2 \in \mathcal{F}$, 
\begin{equation*}
a_1 f_1 + a_2 f_2 = (|a_1| + |a_2|) \left( \frac{a_1}{|a_1| + |a_2|} f_1 +  \frac{a_2}{|a_1| + |a_2|} f_2 \right)     
\end{equation*}
belongs to $\mathcal{F}'$ since by absolute convexity $a_1/(|a_1| + |a_2|) f_1 +  a_2/(|a_1| + |a_2|) f_2 \in \mathcal{F}$. It easily follows that for any $a_1, a_2 \in \R$ and $f_1, f_2 \in \mathcal{F}'$, $a_1 f_1 + a_2 f_2  \in \mathcal{F}'$ and thus $\mathcal{F}'$ is a vector space. As in addition $\mathcal{F} \cup \{ 1 \} \subseteq \mathcal{F}'$, we see that $\mathcal{F}'$ also contains the linear span of $\mathcal{F} \cup \{ 1 \}$. The other inclusion easily holds, and thus there is equality.

\emph{Step 2. Density of the linear span in $C_b(\X)$.} By a corollary of the Hahn-Banach theorem, see Theorem 5.19 in \citet{rudin1987real}, it is enough to prove that any continuous linear form $\sigma$ on $(C_b(\X), \mathrm{L}_\infty)$ that vanishes on $\mathcal{F} \cup \{ 1 \}$ is necessarily identically zero. As $\X$ is compact, the Riesz theorem (see e.g. Theorem 2.14 in \citet{rudin1987real}) guarantees that any such $\sigma$ is characterized by a finite signed measure, which we denote $\sigma$ as well, i.e. $\sigma(f) = \int_\X f \ddr \sigma$. We write $\sigma = \sigma_+ - \sigma_-$ for its Jordan decomposition into a positive and a negative part. As $\sigma(1) = 0$, we deduce that $\sigma_+(\X) = \sigma_-(\X)$. If $\sigma_+(\X) = 0$ then $\sigma_+ = \sigma_- = \sigma = 0$ and we are done. If not, up to dividing $\sigma$ by $\sigma_+(\X)$, we can assume that $\sigma_+(\X) = \sigma_-(\X) = 1$, that is, $\sigma_+$ and $\sigma_-$ are probability distributions over $\X$. Since $\sigma$ vanishes on $\mathcal{F}$, $\sigma_+(f) - \sigma_-(f) = 0$ for any $f \in \mathcal{F}$, which implies $\mathcal{I}_\mathcal{F}(\sigma_+, \sigma_-) = 0$. As $\mathcal{I}_\mathcal{F}$ is distance, we obtain $\sigma_+ = \sigma_-$. This implies $\sigma = 0$ and concludes the proof.

\subsection*{Proof of Proposition~\ref{th:bounds}}

Take $\tilde P_1, \tilde P_2$ distributed according to $\mathbb Q_1$ and $\mathbb Q_2$ respectively. Given the definition of the two distances and the linearity of the integral, the property we want to prove is
\[
\sup_{f \in \mathcal{F}} |\E( \tilde{P}_1(f) ) - \E( \tilde{P}_2(f) )| \le \sup_{h \in \mathfrak F} | \mathbb Q_1(h) - \mathbb Q_2(h) | \le \sup_{h \in \mathrm{Lip}_1(\mathcal{I}_\mathcal{F})} | \mathbb Q_1(h) - \mathbb Q_2(h) |,
\]
where we recall that 
\[
\mathfrak{F} = \{ P \mapsto g(P(f)) \ : \ f \in \mathcal{F} \text{ and } g \in \mathrm{Lip}_1(\R) \}.
\]
Following the definition we see that $\mathfrak{F} \subset \mathrm{Lip}_1(\mathcal{I}_\mathcal{F})$ thus the second inequality holds. Then we notice that for any $f \in \mathcal{F}$ there holds $\E( \tilde{P}_i(f) ) = \mathbb Q_i (h_f)$ for $i=1,2$ with $h_f : P \mapsto P(f) \in \mathfrak{F}$. The first inequality easily follows.

\subsection*{Proof of Lemma~\ref{th:equality}}

Let $P_i = \E( \tilde P_i)$. The fact that $\W(P_1,P_2) \le \W_{\W}( \mathcal{L}(\tilde P_1), \mathcal{L}(\tilde P_2))$ follows from Proposition~\ref{th:bounds} and was also proved in Lemma 3.1 of \citet{Nguyen2016}. We prove that $\W(P_1,P_2) \ge \W_{\W}( \mathcal{L}(\tilde P_1), \mathcal{L}(\tilde P_2))$. Let $T$ be the optimal transport map such that $T_{\#} \tilde P_1 \eqd \tilde P_2$. Then by considering the coupling $(\tilde P_1,T_{\#} \tilde P_1)$, 
\[
 \W_{\W}(\mathcal{L}(\tilde P_1), \mathcal{L}(\tilde P_2)) \le \E_{\tilde P_1} ( \W (\tilde P_1, T_{\#} \tilde P_1 )) = \E_{\tilde P_1} ( \E_{X \mid \tilde P_1 \sim \tilde P_1} (|X - T(X)| \mid \tilde P_1),
 \]
since $T$ is an optimal transport map. If $X|\tilde P_1 \sim \tilde P_1$ then $X \sim \E(\tilde P_1) = P_1$. Thus by the tower property, the right hand side is equal to $\E_{X \sim P_1} (|X - T(X)|)$. We observe that $T_{\#} \tilde P_1 = \tilde P_2$ a.s. implies that $T_{\#}  P_1 =  P_2$. Indeed, $P_2(A) = \E(\tilde P_2(A)) = \E( \tilde P_1( T^{-1}(A)) = P_1(T^{-1}(A)) = T_{\#} P_1(A)$. Since $T$ is an optimal transport map, $\E_{X \sim P_1} (|X - T(X)|) = \W(P_1,P_2)$. Thus, $\W_{\W}(\mathcal{L}(\tilde P_1), \mathcal{L}(\tilde P_2)) \le \W(P_1,P_2)$ and the conclusion follows.

\subsection*{Proof of Theorem~\ref{th:species_sampling}}

If there exists $T$ an optimal transport map between $P_1$ and $P_2$, by assumptions (a), (b), and (c) it follows that
\[
T_{\#} \tilde P_1 = \sum_{j \ge 1} J_j^{(1)} \delta_{T(X_j^{(1)})} \eqd \tilde P_2.
\]
We can apply Lemma~\ref{th:equality} and the conclusion follows.

Consider now the case when an optimal transport map $T$ does not exist and let $\gamma$ be the optimal coupling on $\X \times \X$ between $P_1$ and $P_2$. We build the following coupling between $\law(\tilde P_1)$ and $\law(\tilde P_2)$
\[
(\tilde P_1, \tilde P_2) = \bigg( \sum_{j \ge 1} J_j \delta_{Y_j^{(1)}}, \sum_{j \ge 1} J_j \delta_{Y_j^{(2)}} \bigg),
\]
where $\{J_j\}$ have the same distribution as $\{J_j^{(1)}\}$ and $\{J_j^{(2)}\}$, and they are independent from $(Y_j^{(1)},Y_j^{(2)}) \simiid \gamma$. Then
\begin{align*}
\W_\W( \tilde P_1,\tilde P_2)  & \le \E \bigg( \W \bigg(\sum_{j \ge 1} J_j \delta_{Y_j^{(1)}}, \sum_{j \ge 1} J_j \delta_{Y_j^{(2)}}\bigg) \bigg) \\
& \le  \E \bigg( \sum_{j \ge 1} J_j d_\X(Y_j^{(1)}, Y_j^{(2)}) \bigg) =  \sum_{j \ge 1} \E(J_j) \E(d_\X(Y_j^{(1)}, Y_j^{(2)})).
\end{align*}
Since $(Y_j^{(1)},Y_j^{(2)}) \simiid \gamma$ is the optimal coupling between $P_1$ and $P_2$,
\[
\sum_{j \ge 1} \E(J_j) \E(d(Y_j^{(1)}, Y_j^{(2)})) =  \W(P_1, P_2) \E\bigg( \sum_{j \ge 1}J_j\bigg) = \W(P_1, P_2).
\]
This proves that $\W_\W( \tilde P_1,\tilde P_2) \le \W(P_1, P_2)$. The other inequality follows from Proposition~\ref{th:bounds}.

\subsection*{Proof of Theorem~\ref{th:bad_rates_wow}: the upper bound}

\emph {Step 1: reducing to an estimation of Rademacher complexity.} The quantity of interest can be written, following the definition, 
\[ 
\E(\W_\W(\mathbb{\tilde Q}_{(n)}, \mathbb{Q})) = \E \left( \sup_{h \in \mathrm{Lip}_1^*(\W)} \left| \mathbb{ \tilde Q}_{(n)}(h) - \mathbb Q(h) \right| \right),
\]
where $\rm{Lip}^*_1(\W)$ is the class of functions defined over $\mathcal{P}(\X)$ which are $1$-Lipschitz with respect with the Wasserstein distance and vanish on a distinguished measure $P_0$, which will not play a role. The original definition is stated with $\rm{Lip}_1(\W)$ but clearly replacing by $\rm{Lip}^*_1(\W)$ does not change the value. The gain is that the class $\rm{Lip}^*_1(\W)$ is uniformly bounded as $\X$ is bounded. The classical symmetrization argument (see e.g. Lemma 2.3.1 in \citet{vaartWellner} or Section 4.2 in \citet{Wainwright2019}) yields
\begin{equation*}
\E(\W_\W(\mathbb{\tilde Q}_{(n)}, \mathbb{Q})) \le 2 \mathcal{R}_n( \rm{Lip}^*_1(\W) ),
\end{equation*}
where the Rademacher complexity of the class $\rm{Lip}^*_1(\W)$ is defined as 
\begin{equation}
\label{eq:def_Rademacher}
\mathcal{R}_n(\rm{Lip}^*_1(\W)) = \E_{\tilde P_{1:n}, \epsilon_{1:n}} \bigg( \sup_{h \in \rm{Lip}^*_1(\W)}  \bigg| \frac{1}{n} \sum_{i=1}^n \epsilon_i h(\tilde P_i) \bigg| \bigg),    
\end{equation}
being $\epsilon_1, \ldots, \epsilon_n$ i.i.d. Rademacher random variables independent from $P_1, \ldots, P_n$, which are i.i.d. with law $\mathbb Q$. Note that we use $\E_{\tilde P_{1:n}, \epsilon_{1:n}}$ as a shortcut for $\E_{P_1, \ldots, P_n, \epsilon_1, \ldots, \epsilon_n}$.

\emph{Step 2: estimation of the Rademacher complexity via a covering number.} As trying to get the sharpest rates would take us to expressions which are not analytically tractable, we will only use the one step bound rather than the integral bound (see Proposition 5.17 and Example 5.21 in \citet{Wainwright2019}) which yields
\begin{equation*}
\mathcal{R}_n( \rm{Lip}^*_1(\W) ) \le \inf_{\epsilon \in [0, \epsilon_0)} \left\{ \epsilon + \frac{C_1}{\sqrt{n}} \sqrt{\log N \left( \epsilon, \rm{Lip}^*_1(\W), \rm{L}_\infty  \right)} \right\},    
\end{equation*}
where $\epsilon_0, C_1$ are two constants. Here $N \left( \epsilon, \rm{Lip}^*_1(\W), \rm{L}_\infty  \right)$ is the covering number of the space of $1$-Lipschitz functions defined on $\mathcal{P}(\X)$ with respect to $\rm{L}_\infty$, the supremum norm on the space of functions defined on $\mathcal{P}(\X)$. The original result is stated with $\rm{L}_\infty$ replaced with the $\rm{L}_2$ norm with respect to the empirical process, but is always dominated by $\rm{L}_\infty$. To estimate the covering number we then use the bound of \citet{KolmogorovTikhomirov1961} which is recalled in (3.1) of \citet{sriperumbudur2012}: for a bounded metric space $\mathbb Y$ with metric $d_\Y$, the covering number of $\mathrm{Lip}^*_1(\mathbb Y)$ the one-Lipschitz functions with respect to the uniform norm can be estimated as 
\begin{equation}
\label{eq:bound_KolmogorovTikhomirov}
\log N (\epsilon, \rm{Lip}^*_1(\mathbb{Y}), L_{\infty}) \le N \bigg(\frac{\epsilon}{4}, \mathbb{Y}, d_\Y \bigg)  \log \bigg( 2 \bigg \lceil \frac{2 \rm{diam}(\mathbb{Y})}{\epsilon} \bigg \rceil + 1 \bigg).    
\end{equation}
Here and in the sequel, $\lceil a \rceil$ is the integral value greater than $a$. Considering $(\mathbb Y, d_\Y) = (\mathcal{P}(\X), \W)$ we find
\[
\log N \left( \epsilon, \rm{Lip}^*_1(\W), \rm{L}_\infty  \right) \le N \bigg(\frac{\epsilon}{4}, \mathcal{P}(\X), \W \bigg)  \log \bigg( 2 \bigg \lceil \frac{2 \rm{diam}(\mathcal{P}(\X))}{\epsilon} \bigg \rceil + 1 \bigg).
\]
As $\rm{diam}(\mathcal{P}(\X)) = \rm{diam}(\X)$, what is left to do is to estimate the metric entropy of $\mathcal{P}(\X)$ when endowed with the Wasserstein distance.  

\emph{Step 3: estimating the covering number of the Wasserstein space.} We claim that for any $\epsilon > 0$
\begin{equation}
\label{eq:covering_Wasserstein}
N \bigg(\epsilon, \mathcal{P}(\X), \W \bigg) \le  N \left( \frac{\epsilon}{2},\X,d_\X \right)^{\lceil 2 \rm{diam}(\X) / \epsilon  \rceil +1}.     
\end{equation}
Indeed, let $x_1, x_2, \ldots, x_A$ with $A = N(\epsilon/2,\X,d_\X)$ an $\epsilon/2$-covering of the space $\X$. Let $B = \lceil 2 \rm{diam}(\X) / \epsilon  \rceil$ and let us consider $\mathcal Z$ the subset of $\mathcal{P}(\X)$ made of probability measures supported on $x_1, x_2, \ldots, x_A$ and such that the mass of each $x_j$ belongs to $\{ 0, 1/B, \ldots, (B-1)/B, 1 \}$ for $j=1,\ldots, A$: the cardinality of $\mathcal Z$ is bounded by $A^{B+1}$ and we will prove $\mathcal Z$ is an $\epsilon$-covering of $\mathcal{P}(\X)$.

Take any $P \in \mathcal{P}(\X)$. Then there is a probability measure $P_1$ supported on $x_1, x_2, \ldots, x_A$ such that $\W(P,P_1) \leq \epsilon/2$: it is obtained by projecting the mass of $P$ onto the set $x_1, x_2, \ldots, x_A$. Then, by rounding up the mass of each atom $x_j$, we can find $P_2$ a measure in $\mathcal Z$ such that $\| P_1 - P_2 \|_{\rm{TV}} \le 1/B$, being $\| \cdot \|_{\rm{TV}}$ the total variation norm. Using the easy inequality $\W(P_1, P_2) \le \rm{diam}(\X) \| P_1 - P_2 \|_{\rm{TV}}$ (see Theorem 4 in \citet{gibbs2002choosing}) followed by the triangle inequality, we find that 
\begin{equation*}
\W(P,P_2) \le \W(P,P_1) + \W(P_1,P_2) \le \frac{\epsilon}{2} + \frac{\rm{diam}(\X)}{B} \le \frac{\epsilon}{2} + \frac{\epsilon}{2} = \epsilon.     
\end{equation*}
Thus $\mathcal{Z}$ is an $\epsilon$-covering of $\mathcal{P}(\X)$, and given the expression of its cardinality we obtain the estimate~\eqref{eq:covering_Wasserstein}.

\emph{Step 4: conclusion with a well chosen $\epsilon$}. Chaining our estimates, we see that we can estimate the Rademacher complexity by 
\begin{multline*}
\mathcal{R}_n( \rm{Lip}^*_1(\W) ) \\
 \le \inf_{\epsilon \in [0, \epsilon_0)} \left\{ \epsilon + \frac{C_1}{\sqrt{n}} N \left( \frac{\epsilon}{2},\X,d_\X \right)^{\lceil 2 \rm{diam}(\X) / \epsilon  \rceil /2 +1/2}  \sqrt{\log \bigg( 2 \bigg \lceil \frac{2 \rm{diam}(\mathcal{P}(\X))}{\epsilon} \bigg \rceil + 1 \bigg)} \right\}.    
\end{multline*}
At this point we will not keep track of the explicit constant anymore. As $\X$ is a bounded set of $\R^d$ using Theorem 2.7.1 in \citet{vaartWellner}, we have $N \left( \frac{\epsilon}{2},\X,d_\X \right) \le C_2 \epsilon^{-d}$ for some constant $C_2$. Injecting this in our estimate, for some constant $C_1, C_3, C_4$ large enough and up to decreasing $\epsilon_0$, we get  
\begin{equation}
\label{eq:aux_convergence_rate_wow}
\mathcal{R}_n( \rm{Lip}^*_1(\W) ) \le \inf_{\epsilon \in [0, \epsilon_0)} \left\{ \epsilon + \frac{C_1}{\sqrt{n}} \exp \left( \frac{C_3}{\epsilon} \log \left( \frac{1}{\epsilon} \right)  \right)  \sqrt{\log \left( \frac{C_4}{\epsilon} \right)} \right\}.    
\end{equation}
Exact optimization of this expression in $\epsilon$ seems tricky and leads to intractable analytical expressions. Note however that if $\epsilon \lesssim 1/\log(n)$, then the expression $C_3/\epsilon \log(1/\epsilon)$ in the argument of the exponential grows much faster than $\log(n)$, and so the exponential grows much faster than any positive power of $n$. Thus if $\epsilon \lesssim 1/\log(n)$ the second term in the right hand side of~\eqref{eq:aux_convergence_rate_wow} blows up to $+ \infty$. So we would look for $\epsilon = \epsilon_n = \beta_n / \log(n)$ with $1 \ll \beta_n \ll \log(n)$, and then we see that the argument of the exponential reads 
\begin{equation*}
\frac{C_3}{\epsilon} \log \left( \frac{1}{\epsilon} \right) = C_3 \frac{\log(n)}{\beta_n} \left( \log(\log(n)) - \log(\beta_n) \right) \, \sim C_3 \frac{\log(n) \log(\log(n))}{\beta_n}.     
\end{equation*}
For this to grow slower than $\log(n)$ (so that the exponential grows polynomially) we need $\beta_n$ to grow at least as fast as $\log(\log(n))$. Following this computation we choose
\begin{equation*}
\frac{1}{\epsilon} = \frac{1}{3C_3} \frac{\log(n)}{\log(\log(n))},  
\end{equation*}
which we can do for $n$ large enough, and it will be the best we can do to bound the right hand side of~\eqref{eq:aux_convergence_rate_wow}. We claim that the second term in the right hand side of~\eqref{eq:aux_convergence_rate_wow} is negligible compare to the first one $\epsilon \asymp \log(\log(n))/\log(n)$. Indeed looking at the argument of the exponential, as $1/\epsilon \le \log(n)$ for $n$ large enough, we have 
\begin{equation*}
\frac{C_3}{\epsilon} \log \left( \frac{1}{\epsilon} \right) \le \frac{1}{3} \log(n),     
\end{equation*}
at least for $n$ large enough. Thus the exponential is bounded by $n^{1/3}$ (actually any exponent strictly smaller than $1/2$ would be enough). Moreover $\sqrt{\log(C_4/\epsilon)} \le n^\eta$ for any $\eta > 0$ if $n$ is large enough, and $\eta = 1/12$ will be sufficient for our purposes. So  for $n$ large enough
\begin{equation*}
\frac{C_1}{\sqrt{n}} \exp \left( \frac{C_3}{\epsilon} \log \left( \frac{1}{\epsilon} \right)  \right)  \sqrt{\log \left( \frac{C_4}{\epsilon} \right)}  \le \frac{C_1} {\sqrt{n}} \cdot n^{1/3} \cdot n^{1/12} = C_1 n^{-1/12},    
\end{equation*}
which is negligible compared to the first term which scales like $\log (\log(n))/\log(n)$. Thus
\begin{equation*}
\mathcal{R}_n( \rm{Lip}^*_1(\W) ) \le \frac{C_5 \log (\log(n))}{\log(n)}    
\end{equation*}
for some constant $C_5$, at least if $n$ is large enough. Plugging this into the first step yields the conclusion.

\subsection*{Proof of Theorem~\ref{th:bad_rates_wow}: the lower bound}

Let $P_0$ a measure on $\R^d$ whose support has non-empty interior. Let $\mathbb Q = \text{DP}(\alpha, P_0)$ be a Dirichlet process with parameters $\alpha>0$ and $P_0$. 
For any exponent $\gamma > 0$ we define $D = \lceil 1/\gamma \rceil$. 

The support of $P_0$ contains a ball in $\R^d$. In particular we can find $D$ disjoint and closed sets $A_1,\dots,A_D$ in the support of $P_0$ each of them with non-empty interior. Let $f_1,\dots,f_D: \mathbb{X} \to \mathbb{R}$ be $1$-Lipschitz functions with support in $A_1,\dots,A_D$ respectively and which are not identically zero. E.g., we can set $f_i (\cdot) = d(\cdot,\X \setminus A_i)$ the distance to the complement of $A_i$. We consider the embedding $\iota$ of $\mathcal P(\X)$ into $[0,1]^D$ given by 
\begin{equation*}
\iota(P) = (P(f_1), \ldots, P(f_D)). 
\end{equation*}
Our conclusion will follow by combining the two following claims: (1) as $\iota$ is a Lipschitz map, the rate at which $\E(\W_\W(\tilde{\mathbb Q}_{(n)}, \mathbb Q))$ goes to zero is worse than the statistical rate of convergence for the embedded measures $\iota_\# \tilde{\mathbb Q}_{(n)}$, $\iota_\# \mathbb Q$; (2) the embedded measures live in a Euclidean space of dimension $D$, for which the Wasserstein distance suffers the curse of dimensionality, thus giving a lower bound on the rate.  

\emph {Step 1: analysis of the embedding $\iota$.} The embedding $\iota: \mathcal{P}(\X) \to [0,1]^D$ is a Lipschitz map with respect to the Wasserstein-1 metric and the standard Euclidean distance: this is because the $f_i$ are 1-Lipschitz so that
\begin{align*}
\|\iota(P_1) - \iota(P_2) \| & = \sqrt{ \sum_{i=1}^D | P_1(f_i) - P_2(f_i)|^2 } \\
& \le \sqrt{D} \sup_{f \in \rm{Lip}_1(\X)} | P_1(f) - P_2(f)| = \sqrt{D} \W(P_1,P_2).
\end{align*}
Integrating this inequality with $(\tilde P_1, \tilde P_2)$ any coupling between $(\mathbb{\tilde Q}_{(n)}, \mathbb Q)$, and minimizing over $\mathbb Q$, we find 
\begin{multline*}
\W_{\mathcal{W}}(\mathbb{\tilde Q}_{(n)}, \mathbb{Q}) = \inf_{(\tilde P_{n}, \tilde P) \in \Gamma(\mathbb{\tilde Q}_{(n)}, \mathbb{Q})} \E(\W(\tilde P_{n}, \tilde P)) \\ \ge D^{-1/2} \inf_{(\tilde P_{n}, \tilde P) \in \Gamma(\mathbb{\tilde Q}_{(n)}, \mathbb{Q})} \E( \|i(\tilde P_{n}) - i(\tilde P) \|) \ge D^{-1/2} \W_\mathrm{Euc}(\iota_\# \mathbb{\tilde Q}_{(n)}, \iota_\# \mathbb Q),    
\end{multline*}
where the $\W_\mathrm{Euc}$ denotes the standard $1$-Wasserstein distance on $[0,1]^D$ with the Euclidean distance. Note that $\iota_\# \mathbb{\tilde Q}_{(n)} = (\iota_\# \tilde{\mathbb Q})_{(n)}$ is the empirical estimator of $\iota_\# \mathbb Q$.

\emph{Step 2: analysis of the law of $\iota(\tilde P)$.} Recall that $\tilde{P} \sim \mathbb
Q$ follows a Dirichlet process of parameters $\alpha$ and $P_0$. The law of $\iota(\tilde P) = (\tilde P (f_1), \ldots, \tilde P(f_D))$ as a distribution on $[0,1]^D$ is usually referred as the law of the \emph{vector of random means}. We use Theorem 10 in \citet{lijoi2004means} to prove that $\iota_\# \mathbb Q$, the distribution of $\iota(\tilde P)$, is absolutely continuous with respect to the Lebesgue measure on $\R^D$. This theorem requires $f = (f_1, \ldots f_D)$ to not be affinely $\alpha P_0$-degenerated, that is, for every $(v_1,\dots,v_D) \in \R^D \setminus \{ 0 \}$, the function $x \mapsto \sum_i v_i f_i(x)$ is not $P_0$-a.s. a constant. This assumption is clearly satisfied  as the functions $f_1, \ldots, f_D$ are continuous, not identically zero, with disjoint support, each included in the support of $P_0$. 

\emph{Step 3: conclusion with the curse of dimensionality of the Wasserstein distance}. It is well understood that for an absolutely continuous measure with respect to the Lebesgue measure on $\R^D$ the rate of convergence cannot be better than $n^{-1/D}$ (see, e.g., Proposition 2.1 in \citet{Dudley1969}). Thus, given Step 2, there exists a constant $c_D>0$ such that 
\begin{equation*}
\E \left(\W_\mathrm{Euc}((\iota_\# \tilde{\mathbb Q})_{(n)}, \iota_\# \mathbb Q)  \right) \geq c_D n^{-1/D}.    
\end{equation*}
Combining this with the result of Step 1,
\begin{equation*}
\E \left( \W_{\W}(\mathbb{\tilde Q}_{(n)}, \mathbb{Q}) \right) \geq c_D D^{-1/2} \, n^{-1/D}.    
\end{equation*}
The conclusion follows as $D$ was chosen with $1/D \le \gamma$.

\subsection*{Proof of Lemma~\ref{th:entropy}}

The starting point follows closely the proof of the upper bound of Theorem~\ref{th:bad_rates_wow}, as we reduce to the estimation of a Rademacher complexity and then a covering number. By combining the definition of $d_{\mathcal{F}}$ with the dual formulation of the Wasserstein distance on $\R$ it holds 
\[ 
\E(d_{\mathcal{F}}(\mathbb{\tilde Q}_{(n)}, \mathbb{Q})) = \E \left( \sup_{h \in \mathfrak{F}} \left| \tilde{\mathbb Q}_{(n)}(h) - \mathbb Q(h) \right| \right),
\]
where the class $\mathfrak{F}$ is define by 
\[
\mathfrak{F} = \{ h:\mathcal{P}(\X) \to \mathbb{R} \text{ s.t. } \exists  f \in \mathcal{F}, \exists g \in \rm{Lip}_1(\R,\R) \text{ s.t. } h(P) = g(P(f))\}.
\]

\emph{Step 1: getting back to a bounded class}. To go back to a bounded class we shift the functions by a constant. We define 
\[
\mathfrak{F}^* = \{ h:\mathcal{P}(\X) \to \mathbb{R} \text{ s.t. } \exists  f \in \mathcal{F}^*, \exists g \in \rm{Lip}_1^*([-K,K],\R) \text{ s.t. } h(P) = g(P(f))\},
\]
where $\mathcal{F}^* = \{ f^* = f - f(x_0) \text{ s.t. } f \in \mathcal{F} \}$, with $x_0$ a fixed point in $\X$, $K$ is its uniform bound, i.e., $f(x) \le K$ for every $f \in \mathcal{F}^*$, and $\rm{Lip}_1^*([-K,K],\R) \subseteq \rm{Lip}_1(\R,\R)$ denotes the 1-Lipschitz functions $g:[-K,K] \to \R$ such that $g(0) = 0$. As probability distributions all have the same mass, it is easy to check that the quantity of interest can be expressed as
\[ 
\E(d_{\mathcal{F}}(\mathbb{\tilde Q}_{(n)}, \mathbb{Q})) = \E \left( \sup_{h \in \mathfrak{F}^*} \left| \tilde{\mathbb Q}_{(n)}(h) - \mathbb Q(h) \right| \right). 
\]
Moreover, as the function $g \in \rm{Lip}_1^*([-K,K],\R)$ are obviously bounded by $K$, we see that $\mathcal{F}^*$ uniformly bounded by $K$ over $\X$ translates in $\mathfrak{F}^*$ uniformly bounded by $K$ over $\mathcal{P}(\X)$.   

\emph{Step 2: reducing to an estimation of a Rademacher complexity}. After these preliminary remarks we start similarly to the proof of Theorem~\ref{th:bad_rates_wow}: using the classical symmetrization argument (see e.g. Lemma 2.3.1 in \cite{vaartWellner} or Section 4.2 in \cite{Wainwright2019}) we obtain 
\begin{equation*}
\E(d_{\mathcal{F}}(\mathbb{\tilde Q}_{(n)}, \mathbb{Q})) \le 2 \mathcal{R}_n( \mathfrak{F}^* ),
\end{equation*}
where the Rademacher complexity $\mathcal{R}_n( \mathfrak{F}^* )$ is defined in~\eqref{eq:def_Rademacher}. 

\emph{Step 3: estimating the Rademacher complexity from above with a covering number} By standard arguments using Dudley's entropy integral and the sub-Gaussianity of the Rademacher process (see, e.g., Theorem 5.22 and equation (5.48) in \citet{Wainwright2019}), as the functions in $\mathfrak{F}^*$ are uniformly bounded by $K$ then the Rademacher complexity may be bounded through an integral of the covering number of $\mathfrak{F}^*$ with respect to the empirical $\rm{L}_2$ metric, which is itself controlled by the uniform norm $\rm{L}_\infty$. Let us define $N(\delta; \mathfrak{F}^*, \rm{L}_{\infty})$ denote the $\delta$-covering number of $\mathfrak{F}^*$ with respect to the uniform norm $\rm{L}_{\infty}$. Then the aforementioned results read 
\begin{equation*}
\mathcal{R}_n(\mathfrak{F}^*) \leq \inf_{\epsilon>0} \bigg \{2 \epsilon + \frac{32}{\sqrt{n}} \int_{\epsilon/4}^{K} \sqrt{ \log N(\delta; \mathfrak{F}^*, \rm{L}_{\infty})} \, \ddr \delta \bigg\}.
\end{equation*}
We now bound $N(\delta; \mathfrak{F}^*, \rm{L}_{\infty})$ through the covering number of $\rm{Lip}_1^*([-K,K],\R)$ and $\mathcal{F}^*$ with respect to the supremum norm. Let $g_1,\dots,g_A$ be a $\delta/2$-covering for $\rm{Lip}_1^*([-K,K],\R)$ and let $f_1,\dots,f_B$ be a $\delta/2$-covering for $\mathcal{F}^*$. In particular $A = N(\delta/2,\rm{Lip}_1^*([-K,K],\R),\mathrm{L}_\infty)$ and $B = N(\delta/2,\mathcal{F}^*,\mathrm{L}_\infty)$. We claim that $\{ h_{i,j}(P) = g_i(P(f_j)) \}$ is a $\delta$-covering for $\mathfrak{F}^*$. Indeed, for any $h \in \mathfrak{F}^*$ such that $h(P) = g(P(f))$, let $g_i$ and $f_j$ such that respectively $\| g_i - g \|_{\rm{L}_\infty} \leq \delta/2$ and $\| f_j - f \|_{\rm{L}_\infty} \leq \delta/2$. Then
\begin{align*}
\|h - h_{i,j}\|_{\rm{L}_{\infty}} & \le \sup_P |h(P) - h_{i,j}(P)|\\
& \le \sup_P | g(P(f)) -   g(P(f_j))| + |g(P(f_j)) -  g_i(P(f_j))| \\
& \le \sup_P |P(f) - P(f_j)| + \sup_t | g(t) - g_i(t)| \\
& \le \| f_j - f \|_{\rm{L}_\infty} + \| g_i - g \|_{\rm{L}_\infty} = \frac{\delta}{2} +  \frac{\delta}{2} = \delta.
\end{align*}
This shows that $N(\delta; \mathfrak{F}^*, \rm{L}_{\infty}) \le AB$. The quantity $A$, which is the $\delta/2$ covering number of $\rm{Lip}_1^*([-K,K],\R)$, can be estimated as $\log (A) \le  \log(2) \lceil 4K/\delta \rceil$ (see e.g. p.93 in \citet{KolmogorovTikhomirov1961}). Thus  
\[
\log N(\delta; \mathfrak{F}^*, \rm{L}_{\infty}) \le \log(2) \bigg \lceil \frac{4K}{\delta} \bigg \rceil + \log N \left( \frac{\delta}{2} ; \mathcal{F}^*, \rm{L}_{\infty} \right). 
\]
Plugging this back we obtain 
\[
2 \mathcal{R}_n(\mathfrak{F}^*) \le \inf_{\epsilon>0} \bigg \{4 \epsilon + \frac{64}{\sqrt{n}} \int_{\epsilon/4}^{K} \sqrt{ \log(2) \bigg \lceil \frac{4K}{\delta} \bigg \rceil   + \log N\left( \frac{\delta}{2} ; \mathcal{F}^*, \rm{L}_{\infty} \right)} \, \ddr \delta \bigg\}.
\]
By using the subadditivity of the square root, i.e. $\sqrt{a+b} \le \sqrt{a} + \sqrt{b}$ for all $a,b\ge 0$, the bounds of Step 2 and Step 3 yield
\begin{multline*}
\E(d_{\mathcal{F}}(\mathbb{\tilde Q}_{(n)}, \mathbb{Q})) \\ \le \frac{64 \sqrt{\log(2)}}{\sqrt{n}}  \int_0^K \sqrt{ \bigg \lceil \frac{4K}{\delta}\bigg \rceil } \, \ddr \delta +  \inf_{\epsilon>0} \bigg \{4 \epsilon + \frac{64}{\sqrt{n}} \int_{\epsilon/4}^{K}    \sqrt{ \log N\left( \frac{\delta}{2} ; \mathcal{F}^*, \rm{L}_{\infty} \right)} \, \ddr \delta \bigg\}
\end{multline*}
To finish the proof we need only to estimate the prefactor of the parametric part, which we obtain by computing the first integral on the right hand side. As $\lceil 4K/\delta \rceil \le 4K/\delta + 1$ and again by subadditivity of the square root, it can be estimated from above as
\begin{equation*}
\int_0^K \sqrt{ \bigg \lceil \frac{4K}{\delta}\bigg \rceil } \, \ddr \delta \le \int_0^K \sqrt{ \frac{4K}{\delta} + 1} \le 2 \sqrt{K} \int_0^K \delta^{-1/2} \, \ddr \delta + K = 5K. 
\end{equation*}
The conclusion follows by chaining the estimates.

\subsection*{Proof of Theorem~\ref{th:lipschitz_rates}}

Thanks to Lemma~\ref{th:entropy} the proof now amounts to an evaluation of the covering number of the class $\mathcal{F}^*$, which is standard when studying statistical rates for IPM \cite{Sriperumbudur2010, sriperumbudur2012}. We report the reasoning for completeness.

When $\X$ is a bounded domain of $\R^d$ and $\mathcal{F} = \mathrm{Lip}(\X,\R)$, Theorem 2.7.1 in \citet{vaartWellner} yields the existence of a constant $C_1$ depending on $\rm{diam}(\X)$ such that $\log N (\epsilon, \rm{Lip}_1(\mathbb{X}), L_{\infty}) \le C_1 \epsilon^{-d}$. Substituting in the bound of Lemma~\ref{th:entropy} we obtain that there exists a constant $C_2, C_3$ such that 
\[
\E(\dlip(\mathbb{\tilde Q}_{(n)}, \mathbb{Q})) \le \frac{C_2}{\sqrt{n}} + \inf_{\epsilon>0} \bigg \{ 4 \epsilon + \frac{C_3}{\sqrt{n}} \int_{\epsilon/4}^{2K} \delta^{-d/2} \, \ddr \delta \bigg\}.
\]

To get the announced rate it is enough to optimize the expression in $\epsilon$. For $d = 1$ then $\epsilon = 0$ works as the function $f(\delta) = \delta^{-d/2}$ is integrable in zero; and for $d\ge 2$ we can take $\epsilon = n^{-1/d}$.

\subsection*{Proof of Theorem~\ref{th:hier_estimator}}

By using the triangle inequality,
\[
\E(\dlip(\mathbb{\tilde Q}_{(n,m))}, \mathbb Q)) \le \E(\dlip(\mathbb{\tilde Q}_{(n,m)}, \mathbb{\tilde Q}_{(n)})) + \E(\dlip(\mathbb{\tilde Q}_{(n)}, \mathbb Q )).
\]
Theorem~\ref{th:lipschitz_rates} bounds the second term in the right hand side and so we only have to show that for every $n$,
\[
\E(d_{\mathcal{F}}(\mathbb{\tilde Q}_{(n,m)}, \mathbb{\tilde Q}_{(n)})) \le 
\begin{cases}
C_1 m^{-1/2} \qquad &\text{if } d = 1,\\
C_2 m^{-1/2} \log(m) \qquad &\text{if } d = 2,\\
C_d m^{-1/d} \qquad &\text{if }  d >2,
\end{cases} 
\]
where $C_d$ does not depend on $n$ and $m$. By definition of HIPM,
\[
\dlip(\mathbb{\tilde Q}_{(n,m)}, \mathbb{\tilde Q}_{(n)}) = \sup_{f \in \mathcal{F}} \W \bigg( \frac{1}{n} \sum_{i=1}^n \delta_{\tilde P_{i,(m)}(f)}, \frac{1}{n} \sum_{i=1}^n \delta_{\tilde P_{i}(f)} \bigg). 
\]
By coupling $\tilde P_{i,(m)}(f)$ with $\tilde P_{i}(f)$ we obtain a natural upper bound 
\begin{align*}
\E( \dlip(\mathbb{\tilde Q}_{(n,m)}, \mathbb{\tilde Q}_{(n)})) & \le \E \bigg( \sup_{f \in \mathcal{F}} \frac{1}{n} \sum_{i=1}^n |\tilde P_{i,(m)}(f)-\tilde P_{i}(f)|\bigg) \\
& = \frac{1}{n} \sum_{i=1}^n \E \bigg( \sup_{f \in \mathcal{F}}|\tilde P_{i,(m)}(f)-\tilde P_{i}(f)|\bigg).
\end{align*}
By Corollary~8 in \citet{Sriperumbudur2010} on the estimation of the Wasserstein distance of order 1 through the empirical distribution, 
\[
\E \bigg( \sup_{f \in \mathcal{F}}|\tilde P_{i,(m)}(f)-\tilde P_{i}(f)| \bigg|  \tilde P_{i} \bigg) \le 
\begin{cases}
C_1 m^{-1/2} \qquad &\text{if } d = 1,\\
C_2 m^{-1/2} \log(m) \qquad &\text{if } d = 2,\\
C_d m^{-1/d} \qquad &\text{if }  d >2,
\end{cases} 
\]
whenever $\X \subset \R^d$ is a bounded convex set with non-empty interior. The same inequality holds for any bounded set in $\R^d$ as it can be embedded in a bounded convex set with non-empty interior. Since the upper bound does not depend on $\tilde P_i$, by the towering property we obtain the desired upper bound.

\subsection*{Proof of Proposition~\ref{th:approx}}

\emph{Upper bounds for the Dirichlet multinomial $\tilde P_1$.} Note that we can write the law of $\tilde P_1$ as 
\begin{equation*}
\law (\tilde P_1) = \E_{X_{1:N}} \left( \rm{DP}\bigg(\alpha, \frac{1}{N} \sum_{i=1}^N \delta_{X_i}\bigg) \right),    
\end{equation*}
where we use $X_{1:N}$ as shortcut for $X_1, \ldots, X_N \simiid P_0$.
By convexity of the Wasserstein distance (Theorem 4.8 in \citet{Villani2008}),
\begin{align*}
\W_\W(\tilde P, \tilde P_1) & \le  \E_{X_{1:N}} \bigg( \W_\W \bigg( \rm{DP}(\alpha, P_0),  \rm{DP}\bigg(\alpha, \frac{1}{N} \sum_{i=1}^N \delta_{X_i}\bigg) \bigg) \\
& = \E_{X_{1:N}} \bigg(\W \bigg( P_0, \frac{1}{N} \sum_{i=1}^N \delta_{X_i} \bigg)\bigg),
\end{align*}
where the last equality holds by Corollary~\ref{th:exact_dp}. The convergence rate of this quantity is well-studied in the literature when $\X = \R^d$ (see, e.g., \citet{Dudley1969, BoissardLeGouic2014, FournierGuillin2015, WeedBach2019, BobkovLedoux2019} and references therein). In particular, if $\X = \R$ Theorem 3.2 in \citet{BobkovLedoux2019} guarantees that 
\[
\E\bigg(\W \bigg( P_0, \frac{1}{N} \sum_{i=1}^N \delta_{X_i} \bigg)\bigg) \le \frac{1}{\sqrt{N}} \int_{-\infty}^{+\infty} \sqrt{F_0(x)(1- F_0(x))} \ddr x.
\]

\emph{Upper bounds for the truncated stick breaking $\tilde P_2$}. Given a random sequence $(J_i)_{i \geq 1}$ of stick-breaking weights, and a sequence $(X_i)_{i \geq 1}$ of i.i.d. atoms, we can naturally build a coupling between $\mathcal{L}(\tilde P)$ and $\mathcal{L}(\tilde P_2)$ via
\begin{equation*}
\left( \sum_{i\ge 1} J_i \delta_{X_i}, \, \sum_{i=1}^{N-1} J_i \delta_{X_i} + \bigg(1 - \sum_{i=1}^{N-1} J_i \bigg) \delta_{X_N} \right).    
\end{equation*}
This provides a natural upper bound for the Wasserstein over Wasserstein distance. Using $\sum_i J_i = 1$ a.s., we rewrite it as
\[
\W_\W(\tilde P, \tilde P_2) \le \E \bigg( \W \bigg( \sum_{i=1}^{N-1} J_i \delta_{X_i} + \sum_{i=N}^{+\infty} J_i \delta_{X_i}, \sum_{i=1}^{N-1} J_i \delta_{X_i} + \bigg(\sum_{i=N}^{+\infty} J_i  \bigg) \delta_{X_N}\bigg) \bigg). 
\]
By using the representation of the Wasserstein distance as an IPM, 
\begin{multline*}
\W \bigg( \sum_{i=1}^{N-1} J_i \delta_{X_i} + \sum_{i=N}^{+\infty} J_i \delta_{X_i}, \sum_{i=1}^{N-1} J_i \delta_{X_i} + \bigg(\sum_{i=N}^{+\infty} J_i  \bigg) \delta_{X_N}\bigg) \\
 = \sup_{f \in \mathrm{Lip}_1(\X)} \left| \sum_{i=N+1}^{+ \infty} J_i (f(X_i) - f(X_N)) \right|. 
\end{multline*}
We then use the natural bound $|f(X_i) - f(X_N)| \leq d_\X(X_i,X_N)$ independent of $f$. Note that $(X_i, X_N)$ has the same law as $(X_1,X_2)$ by the i.i.d. assumption. Thus when we take the expectation:
\begin{align*}
\W_\W(\tilde P, \tilde P_2) & \leq \E \left( \sum_{i=N+1}^{+ \infty} J_i d_\X(X_i,X_N)  \right) \\
& = \E(d_\X(X_1,X_2)) \E \left( \sum_{i=N+1}^{+ \infty} J_i \right) = \E(d_\X(X_1,X_2)) \left( \frac{\alpha}{\alpha+1} \right)^N,
\end{align*}
where the last equality follows from Section 3.2 of \citet{IshwaranJames2001}. Moreover, when $\X = \R$ the following identity holds (cfr., e.g., \citet{BobkovLedoux2019}) 
\[
\E(|X_1 - X_2|)  = 2 \int_{-\infty}^{+\infty} F_0(x)(1- F_0(x)) \ddr x.
\]

\emph{Upper bounds for the hierarchical empirical measure $\tilde P_3$}. The coupling that sends $\tilde P$ to the Dirichlet process that defines the law of the atoms of $\tilde P_3$ ensures that
\[
\W_\W(\tilde P, \tilde P_3) \le \E_{\tilde P} \bigg( \W \bigg(\tilde P, \frac{1}{N}\sum_{i=1}^N \delta_{X_i} \bigg) \bigg) = \E_{\tilde P} \bigg( \E \bigg( \W \bigg( \tilde P, \frac{1}{N}\sum_{i=1}^N \delta_{X_i} \bigg) \bigg| \tilde P\bigg) \bigg), 
\]
thanks to the tower property. By Theorem 3.2 in \citet{BobkovLedoux2019} it follows that 
\begin{align*}
\W_\W(\tilde P, \tilde P_3) & \le \E_{\tilde P} \bigg( \frac{1}{\sqrt{N}} \int_{-\infty}^{+\infty} \sqrt{\tilde F(x)(1- \tilde F(x))} \ddr x \bigg)
\\ 
& \le  \frac{1}{\sqrt{N}} \int_{-\infty}^{+\infty} \sqrt{ \E( \tilde F(x)) - \E(\tilde F(x)^2)} \ddr x,
\end{align*}
where the last inequality holds by linearity of the expectation, Fubini's Theorem, and Jensen's inequality. Standard properties of the Dirichlet process ensure that $\E( \tilde F(x))=F_0(x)$ and $\E(\tilde F(x)^2) = \text{Var}(\tilde F(x)) + \E(F(x))^2 = (F_0(x)(1-F_0(x)))/(1+\alpha) + F_0(x)^2$. Thus the last term is equal to
\[
\sqrt{\frac{\alpha}{N(\alpha+1)}} \int_{-\infty}^{+\infty} \sqrt{F_0(x)(1- F_0(x))} \ddr x.
\]

\section{Additional information on the algorithms}
\label{app:algorithm}

\subsection{Wasserstein over Wasserstein}

In Algorithm~\ref{alg:WoW} we present the pseudocode for computing the Wasserstein over Wasserstein distance as described in Section~\ref{sec:numerics}.

\begin{algorithm}
   \caption{Computation of $\W_\W$}
   \label{alg:WoW}
\begin{algorithmic}
    \STATE {\bfseries Routine needed}: \texttt{OT}($A$,$a$,$b$) which outputs the value of the transport problem with cost matrix $A$ of size $q \times q$ and weights $a$, $b$ for the marginals. 
   \STATE {\bfseries Input:} sizes $n,m$, data $X^{1}_{i,j}$, $X^{2}_{i,j}$, $\omega^{1}_{i,j}$, $\omega^{2}_{i,j}$
    \STATE Initialize $C=0$ of size $n \times n$
    \algorithmiccomment{Compute the cost matrix for the $n \times n$ transport problem}
    \FOR{$i_1,i_2=1$ {\bfseries to} $n$}
        \STATE $B \leftarrow $ \texttt{PairwiseDistance}($X^{1}_{i_1,\cdot}$, $X^{2}_{i_2,\cdot}$)
        \STATE $C_{i_1, i_2} \leftarrow $ \texttt{OT}($B$, $\omega^{1}_{i_1,\cdot}$, $\omega^{2}_{i_2,\cdot}$) 
    \ENDFOR

    \STATE {\bfseries Return:} \texttt{OT}($C$, $1/n \, \mathbf{1}$, $1/n \, \mathbf{1}$)
    \algorithmiccomment{Solve the $n \times n$ transport problem}
\end{algorithmic}
\end{algorithm}

\subsection{Lipschitz HIPM in dimension one}

We expand on the gradient ascent algorithm for our new distance $\dlip$. Recall that we are solving the optimization problem: 
\begin{align}
\label{eq:problem_to_optimize}
\sup_{\mathbf{f} \in \R^{\ngr}} & \; \mathcal{G}(\mathbf{f})  \\
\notag
\text{such that} & \; |\mathbf{f}_{q+1} - \mathbf{f}_{q}| \leq \Delta x \; \forall q \in \{ 1, \ldots, {\ngr}-1 \}, \\
\notag
\text{where} & \; \mathcal{G}(\mathbf{f}) :=  \inf_{\sigma \in \mathcal{S}(n)} \frac{1}{n} \sum_{i=1}^n \left| \sum_{q=1}^{\ngr} (\omega^1_{i,q} - \omega^2_{\sigma(i),q}) \mathbf{f}_q  \right|.
\end{align}
The function $\mathcal{G}$ is piece-wise linear. On each ``facet'', that is, when there exists a unique permutation $\sigma^*$ which realizes the infimum, and when each term $(\omega^1_{i,q} - \omega^2_{\sigma(i),q}) \mathbf{f}_q$ is either strictly positive or strictly negative, then we can easily compute the gradient. It reads: for any coordinate $q$,
\begin{equation*}
\nabla_q \mathcal{G}(\mathbf{f}) = 
\frac{1}{n} \sum_{i=1}^n (\omega^1_{i,q} - \omega^2_{\sigma^*(i),q}) \, \mathrm{sign}\left( \sum_{q'=1}^N (\omega^1_{i,q'} - \omega^2_{\sigma^*(i),q'}) \mathbf{f}_q \right) .
\end{equation*}
If $\sigma^*$ is not unique or $(\omega^1_{i,q} - \omega^2_{\sigma(i),q})$ vanishes 
the gradient may not exist, but in practice we ignore these degenerate cases. \\
We further rewrite the problem by parametrizing the function $\mathbf{f}$ rather by its derivative to simplify the Lipschitz constraint. Let $\mathbf{g} \in \R^{\ngr-1}$ and consider the $\ngr \times (\ngr -1)$ matrix 
\begin{equation*}
A = \Delta x \begin{pmatrix}
0 & 0 & \ldots & \ldots & 0   \\
1 & 0 & \ldots & \ldots & 0 \\
1 & 1 & 0 & \ddots & \vdots \\
\vdots & & \ddots & & 0 \\
 1 & 1 & \ldots & 1 & 1
\end{pmatrix}.
\end{equation*}

We consider $\mathbf f = A \mathbf g$, that is, $\mathbf{f}_q = \sum_{q' < q} \Delta x \, \mathbf{g}_{q'}$ for any $q$. Specifically we define $\hat{\mathcal{G}}(\mathbf{g}) = \mathcal{G}(A \mathbf f)$, so that problem~\eqref{eq:problem_to_optimize} can be rewritten as 
$\sup \hat{\mathcal{G}}( \mathbf{g})$ for $\mathbf{g} \in [-1,1]^{\ngr-1}$. The gain is that the constraint on $\mathbf{g}$ is a very simple box constraint. Moreover the gradient can be easily computed via $\nabla \hat{\mathcal{G}} (\mathbf g) = A^\top \nabla \mathcal{G}(A \mathbf f)$.

We implemented the following gradient ascent. Given $\mathbf{g}$ admissible we find an ascent direction $\mathbf{a}$ by orthogonally projecting $\nabla \hat{\mathcal{G}}(\mathbf{g})$ on the set of vectors that preserves the box constraint. This can be easily done: starting with $\mathbf{a} = \nabla \hat{\mathcal{G}}(\mathbf{g})$, we set $\mathbf a_q = 0$ if $\mathbf{g}_q = 1$ and $\nabla_q \hat{\mathcal{G}}(\mathbf{g}) > 0$ (respectively if $\mathbf{g}_q = -1$ and $\nabla_q \hat{\mathcal{G}}(\mathbf{g}) < 0$). We then find $t_\mathrm{max}$ the largest $t \geq 0$ such that $\mathbf{g} + t \mathbf{a} \in [-1,1]^{\ngr}$, and we use a backtracking linesearch to find $t \in [0,t_\mathrm{max}]$ such that $\mathcal{G}(\mathbf{g} + t \mathbf{a})$ increases enough \citep{armijo1966minimization}. We stop the loop if the expected increase, that is, $\mathbf{a}^\top \nabla \hat{\mathcal{G}}(\mathbf{g})$, is too small. 

As we have no guarantee of finding a global maximizer we run the gradient ascents with several different initializations. We usually include the initialization $f = \mathrm{Id}$, that is, $\mathbf g = (1,1,\ldots, 1)$ in at least one of the runs, which is related to the bound~\eqref{eq:lower}. This is summarized in the Algorithm~\ref{alg:df}.

\begin{algorithm}
   \caption{Computation of an approximation of $\dlip$}
   \label{alg:df}
\begin{algorithmic}
    \STATE {\bfseries Input:} sizes $n,\ngr$, data $\omega^{1}_{i,q}$, $\omega^{2}_{i,q}$, stepsize $\Delta x$, 
    \STATE {\bfseries Parameters:} number of initializations $n_\text{init}$, number of steps $n_\text{step}$, tolerance $\varepsilon$
    \FOR{$s=1$ {\bfseries to} $n_\text{init}$}
        \STATE $\mathbf{g} \leftarrow $ \texttt{Random Initialization} or $\mathbf{g} \leftarrow (1,1,\ldots,1)$
        \REPEAT 
            \STATE $\mathbf{a} \leftarrow \nabla \hat{\mathcal{G}}(\mathbf{g})$ 
            \algorithmiccomment{Ascent direction}
            \FOR{$q=1$ {\bfseries to} $\ngr-1$}
            \STATE {\bfseries if } $\mathbf{g}_q = \pm 1$ {\bfseries then } $\pm \mathbf{a}_q = \min(0, \pm \mathbf{a}_q) $ {\bfseries end if} 
            \algorithmiccomment{Projection of the ascent direction}
            \ENDFOR 
            \STATE $t_\mathrm{max} \leftarrow \sup \{ t \geq 0 \ : \ \mathbf{g} + t \mathbf{a} \} \in [-1,1]^{\ngr-1}$
            \STATE $t \leftarrow t_\mathrm{max}$ 
            \algorithmiccomment{Backtracking line search}
            \WHILE{$\mathcal{G}(\mathbf g + t \mathbf{a}) < \mathcal{G}(\mathbf{g}) +  t \mathbf{a}^\top \nabla \hat{\mathcal{G}}(\mathbf{g}) /2 $} 
                \STATE $t \leftarrow t/2$ 
            \ENDWHILE
            \STATE $\mathbf{g} \leftarrow \mathbf{g} + t \mathbf{a}$
            \algorithmiccomment{Update of the gradient ascent}
        \UNTIL{ $\mathbf{a}^\top \nabla \hat{\mathcal{G}}(\mathbf{g}) \leq \varepsilon$  {\bfseries or} loop done $n_\text{step}$ times }
        \STATE Store $ \hat{\mathcal{G}}(\mathbf{g})$
    \ENDFOR 
    \STATE {\bfseries Return:} Maximum $ \hat{\mathcal{G}}(\mathbf{g})$ among the $n_\text{init}$ runs 
\end{algorithmic}
\end{algorithm}

\subsection{Execution time}

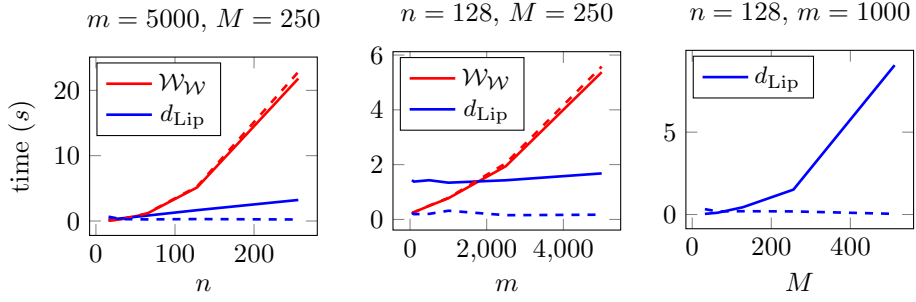
\begin{figure}
\begin{center}
\begin{tabular}{ccc}

\begin{tikzpicture}

\begin{axis}[
width=0.35\textwidth,
xlabel={$n$},
ylabel={time ($s$)},
title ={$m=5000$, $\ngr = 250$},
legend pos= north west]

\addplot[color=red, line width = 1pt] table [x=nTop, y=splitWoW]{timing_nTop.txt};
\addlegendentry{{\small $\W_\W$}}

\addplot[color=blue, line width = 1pt] table [x=nTop, y=splitND]{timing_nTop.txt};
\addlegendentry{{\small $\dlip$}}

\addplot[color=red, line width = 1pt, style = dashed] table [x=nTop, y=sameWoW]{timing_nTop.txt};

\addplot[color=blue, line width = 1pt, style = dashed] table [x=nTop, y=sameND]{timing_nTop.txt};

\end{axis}
\end{tikzpicture}

& 

\begin{tikzpicture}
\begin{axis}[
width=0.35\textwidth,
xlabel={$m$},
title ={$n=128$, $\ngr = 250$},
legend pos= north west]

\addplot[color=red, line width = 1pt] table [x=nBottom, y=splitWoW]{timing_nBottom.txt};
\addlegendentry{{\small $\W_\W$}}

\addplot[color=blue, line width = 1pt] table [x=nBottom, y=splitND]{timing_nBottom.txt};
\addlegendentry{{\small $\dlip$}}

\addplot[color=red, line width = 1pt, style = dashed] table [x=nBottom, y=sameWoW]{timing_nBottom.txt};

\addplot[color=blue, line width = 1pt, style = dashed] table [x=nBottom, y=sameND]{timing_nBottom.txt};

\end{axis}
\end{tikzpicture}

&

\begin{tikzpicture}
\begin{axis}[
width=0.35\textwidth,
xlabel={$\ngr$},
title ={$n=128$, $m = 1000$},
legend pos= north west]

\addplot[color=blue, line width = 1pt] table [x=nGrid, y=splitND]{timing_nGrid.txt};
\addlegendentry{{\small $\dlip$}}

\addplot[color=blue, line width = 1pt, style = dashed] table [x=nGrid, y=sameND]{timing_nGrid.txt};

\end{axis}
\end{tikzpicture}

\end{tabular}
\end{center}
\caption{Execution time (in seconds) for several configurations of the parameters $n$, $m$ and $\ngr$. The setting corresponds to the one of Figure~\ref{fig:distance_numerics} left (for the solid lines) and Figure~\ref{fig:distance_numerics} right (for the dashed lines). Computations are repeated 12 times and the averaged execution time is reported.}
\label{fig:execution_time}
\end{figure}

Computations were performed on the CPU of a standard laptop with a 6-core 2.10GHz AMD Ryzen 5 5500U processor with 8Go of RAM. We report in Figure~\ref{fig:execution_time} the execution time of the computation of $\W_\W$ and (the approximation of) $\dlip$. We do so for different values of $n$, $m$ and also $\ngr$ the number of grid points in the setting of Figure~\ref{fig:distance_numerics}. Whereas the time for computing $\W_\W$ grows quickly with $n$ and is quite insensible to the input measures, we see on the other hand that $\ngr$ and the input measures are the sensible parameters for $\dlip$. The input measures determine the optimal $\mathbf{f}$, and this is likely to affect the convergence of the gradient ascent. Though $\dlip$ is unsurprisingly slow to compute when $\ngr$ is large, recall from that from our error analysis it is natural to take $\ngr \asymp m^{1/2}$, thus moderately large.

\bibliography{references}

\begin{thebibliography}{51}
\providecommand{\natexlab}[1]{#1}
\providecommand{\url}[1]{{#1}}
\providecommand{\urlprefix}{URL }
\providecommand{\doi}[1]{\url{https://doi.org/#1}}
\providecommand{\eprint}[2][]{\url{#2}}
 \bibcommenthead

\bibitem[{Alvarez-Melis and Fusi(2020)}]{AlvarezMelis2020}
Alvarez-Melis D, Fusi N (2020) Geometric dataset distances via optimal
  transport. In: Larochelle H, Ranzato M, Hadsell R, et~al (eds) Advances in
  Neural Information Processing Systems, vol~33. Curran Associates, Inc., pp
  21428--21439

\bibitem[{Ambrosio et~al(2008)Ambrosio, Gigli, and
  Savar{\'e}}]{ambrosio2005gradient}
Ambrosio L, Gigli N, Savar{\'e} G (2008) Gradient flows: in metric spaces and
  in the space of probability measures. Springer Science \& Business Media

\bibitem[{Arbel et~al(2019)Arbel, De~Blasi, and Pr{\"u}nster}]{Arbel2019}
Arbel J, De~Blasi P, Pr{\"u}nster I (2019) {Stochastic Approximations to the
  Pitman–Yor Process}. Bayesian Analysis 14(4):1201 -- 1219

\bibitem[{Armijo(1966)}]{armijo1966minimization}
Armijo L (1966) {Minimization of functions having Lipschitz continuous first
  partial derivatives}. Pacific Journal of mathematics 16(1):1--3

\bibitem[{Bing et~al(2022{\natexlab{a}})Bing, Bunea, and
  Niles-Weed}]{bing2022sketched}
Bing X, Bunea F, Niles-Weed J (2022{\natexlab{a}}) {Estimation and inference
  for the Wasserstein distance between mixing measures in topic models}. arXiv
  preprint arXiv:220612768

\bibitem[{Bing et~al(2022{\natexlab{b}})Bing, Bunea, Strimas-Mackey, and
  Wegkamp}]{bing2022annals}
Bing X, Bunea F, Strimas-Mackey S, et~al (2022{\natexlab{b}}) {Likelihood
  estimation of sparse topic distributions in topic models and its applications
  to Wasserstein document distance calculations}. The Annals of Statistics
  50(6):3307 -- 3333

\bibitem[{Blackwell and MacQueen(1973)}]{Blackwell1973}
Blackwell D, MacQueen JB (1973) {Ferguson distributions via P{\'o}lya urn
  schemes}. The Annals of Statistics 1(2):353--355

\bibitem[{Bobkov and Ledoux(2019)}]{BobkovLedoux2019}
Bobkov SG, Ledoux M (2019) One-dimensional empirical measures, order
  statistics, and {K}antorovich transport distances. Memoirs of the American
  Mathematical Society

\bibitem[{Boissard and Gouic(2014)}]{BoissardLeGouic2014}
Boissard E, Gouic TL (2014) {On the mean speed of convergence of empirical and
  occupation measures in Wasserstein distance}. Annales de l'Institut Henri
  Poincaré, Probabilités et Statistiques 50(2):539 -- 563

\bibitem[{Campbell et~al(2019)Campbell, Huggins, How, and
  Broderick}]{Campbell2019}
Campbell T, Huggins JH, How JP, et~al (2019) {Truncated random measures}.
  Bernoulli 25(2):1256 -- 1288

\bibitem[{Deshpande et~al(2019)Deshpande, Hu, Sun, Pyrros, Siddiqui, Koyejo,
  Zhao, Forsyth, and Schwing}]{deshpande2019max}
Deshpande I, Hu YT, Sun R, et~al (2019) {Max-sliced Wasserstein Distance and
  its use for GANs}. In: Proceedings of the IEEE/CVF conference on computer
  vision and pattern recognition, pp 10648--10656

\bibitem[{Dudley(1969)}]{Dudley1969}
Dudley RM (1969) {The speed of mean Glivenko-Cantelli Convergence}. The Annals
  of Mathematical Statistics 40(1):40 -- 50

\bibitem[{Dukler et~al(2019)Dukler, Li, Lin, and Montufar}]{Dukler2019}
Dukler Y, Li W, Lin A, et~al (2019) {W}asserstein of {W}asserstein loss for
  learning generative models. In: Chaudhuri K, Salakhutdinov R (eds)
  Proceedings of the 36th International Conference on Machine Learning,
  Proceedings of Machine Learning Research, vol~97. PMLR, pp 1716--1725

\bibitem[{Ferguson(1973)}]{Ferguson1973}
Ferguson TS (1973) A {B}ayesian analysis of some nonparametric problems. The
  Annals of Statistics 1(2):209 -- 230

\bibitem[{Fournier and Guillin(2015)}]{FournierGuillin2015}
Fournier N, Guillin A (2015) On the rate of convergence in {W}asserstein
  distance of the empirical measure. Probability Theory and Related Fields
  162:707–738

\bibitem[{Ghosal and van~der Vaart(2017)}]{Ghosal2017}
Ghosal S, van~der Vaart A (2017) Fundamentals of nonparametric Bayesian
  inference. Cambridge University Press

\bibitem[{Gibbs and Su(2002)}]{gibbs2002choosing}
Gibbs AL, Su FE (2002) On choosing and bounding probability metrics.
  International statistical review 70(3):419--435

\bibitem[{Harris(1971)}]{Harris1971}
Harris T (1971) Random measures and motions of point processes. Zeitschrift
  f{\"u}r Wahrscheinlichkeitstheorie und verwandte Gebiete 18:85–115

\bibitem[{Ho et~al(2017)Ho, Nguyen, Yurochkin, Bui, Huynh, and Phung}]{Ho2017}
Ho N, Nguyen X, Yurochkin M, et~al (2017) Multilevel clustering via
  {W}asserstein means. In: International conference on machine learning, PMLR,
  pp 1501--1509

\bibitem[{Ishwaran and James(2001)}]{IshwaranJames2001}
Ishwaran H, James LF (2001) Gibbs sampling methods for stick-breaking priors.
  Journal of the American Statistical Association 96(453):161--173

\bibitem[{Ishwaran and Zarepour(2000)}]{IshwaranZarepour2000}
Ishwaran H, Zarepour M (2000) Markov chain monte carlo in approximate dirichlet
  and beta two-parameter process hierarchical models. Biometrika 87(2):371--390

\bibitem[{Ishwaran and Zarepour(2002)}]{IshwaranZarepour2002}
Ishwaran H, Zarepour M (2002) Exact and approximate sum representations for the
  dirichlet process. The Canadian Journal of Statistics / La Revue Canadienne
  de Statistique 30(2):269--283

\bibitem[{Kallenberg(2017)}]{Kallenberg2017}
Kallenberg O (2017) Random Measures, Theory and Applications. Probability
  Theory and Stochastic Modelling, Springer

\bibitem[{Kolmogorov and Tikhomirov(1961)}]{KolmogorovTikhomirov1961}
Kolmogorov AN, Tikhomirov VM (1961) $\epsilon$-entropy and $\epsilon$-capacity
  of sets in functional spaces. American Mathematical Society Translations
  2(17):277–364

\bibitem[{Lijoi and Regazzini(2004)}]{lijoi2004means}
Lijoi A, Regazzini E (2004) {Means of a Dirichlet process and multiple
  hypergeometric functions}. The Annals of Probability 32(2):1469--1495

\bibitem[{Lijoi et~al(2020)Lijoi, Prünster, and Rigon}]{Lijoi2020}
Lijoi A, Prünster I, Rigon T (2020) {The Pitman–Yor multinomial process for
  mixture modelling}. Biometrika 107(4):891--906

\bibitem[{Mena and Niles-Weed(2019)}]{mena2019statistical}
Mena G, Niles-Weed J (2019) {Statistical bounds for entropic optimal transport:
  sample complexity and the central limit theorem}. Advances in neural
  information processing systems 32

\bibitem[{Müller(1997)}]{Mueller1997}
Müller A (1997) Integral probability metrics and their generating classes of
  functions. Advances in Applied Probability 29(2):429--443

\bibitem[{Muliere and Tardella(1998)}]{MuliereTardella1998}
Muliere P, Tardella L (1998) Approximating distributions of random functionals
  of {F}erguson-{D}irichlet priors. The Canadian Journal of Statistics / La
  Revue Canadienne de Statistique 26(2):283--297

\bibitem[{Nguyen et~al(2023)Nguyen, Huggins, Masoero, Mackey, and
  Broderick}]{Nguyen2023}
Nguyen TD, Huggins J, Masoero L, et~al (2023) Independent finite approximations
  for {B}ayesian nonparametric inference. Bayesian Analysis pp 1 -- 38

\bibitem[{Nguyen(2016)}]{Nguyen2016}
Nguyen X (2016) {Borrowing strengh in hierarchical Bayes: Posterior
  concentration of the Dirichlet base measure}. Bernoulli 22(3):1535 -- 1571

\bibitem[{Orlin(1997)}]{orlin1997polynomial}
Orlin JB (1997) A polynomial time primal network simplex algorithm for minimum
  cost flows. Mathematical Programming 78:109--129

\bibitem[{Peyr{\'e} and Cuturi(2019)}]{peyre2019computational}
Peyr{\'e} G, Cuturi M (2019) {Computational optimal transport: With
  applications to data science}. Foundations and Trends{\textregistered} in
  Machine Learning 11(5-6):355--607

\bibitem[{Pitman(1996)}]{Pitman1996}
Pitman J (1996) Some developments of the {B}lackwell-{M}acqueen urn scheme.
  Lecture Notes-Monograph Series 30:245--267

\bibitem[{Pitman and Yor(1997)}]{PitmanYor1997}
Pitman J, Yor M (1997) {The two-parameter Poisson-Dirichlet distribution
  derived from a stable subordinator}. The Annals of Probability 25(2):855 --
  900

\bibitem[{Prohorov(1961)}]{Prohorov1961}
Prohorov Y (1961) Random measures on a compactum. Soviet Mathematics Doklady
  2:539--541

\bibitem[{Regazzini et~al(2003)Regazzini, Lijoi, and
  Pr\"{u}nster}]{Regazzini2003}
Regazzini E, Lijoi A, Pr\"{u}nster I (2003) Distributional results for means of
  normalized random measures with independent increments. The Annals of
  Statistics 31(2):560--585

\bibitem[{Rudin(1987)}]{rudin1987real}
Rudin W (1987) {Real and Complex Analysis}. Mathematics series, McGraw-Hill

\bibitem[{Sethuraman(1994)}]{Sethuraman1994}
Sethuraman J (1994) {A constructive definition of Dirichlet priors}. Statistica
  Sinica 4(2):639 -- 650

\bibitem[{Sriperumbudur et~al(2010)Sriperumbudur, Fukumizu, Gretton,
  Sch{\"o}lkopf, and Lanckriet}]{Sriperumbudur2010}
Sriperumbudur BK, Fukumizu K, Gretton A, et~al (2010) Non-parametric estimation
  of integral probability metrics. In: 2010 IEEE International Symposium on
  Information Theory, pp 1428--1432

\bibitem[{Sriperumbudur et~al(2011)Sriperumbudur, Fukumizu, and
  Lanckriet}]{sriperumbudur2011universality}
Sriperumbudur BK, Fukumizu K, Lanckriet GR (2011) {Universality, Characteristic
  Kernels and RKHS Embedding of Measures}. Journal of Machine Learning Research
  12(7)

\bibitem[{Sriperumbudur et~al(2012)Sriperumbudur, Fukumizu, Gretton,
  Sch{\"o}lkopf, and Lanckriet}]{sriperumbudur2012}
Sriperumbudur BK, Fukumizu K, Gretton A, et~al (2012) On the empirical
  estimation of integral probability metrics. Electronic Journal of Statistics
  6:550–1599

\bibitem[{Vaart and Wellner(1996)}]{vaartWellner}
Vaart Avd, Wellner JA (1996) Weak Convergence and Empirical Processes: With
  Applications to Statistics. Springer

\bibitem[{Varadarajan(1958)}]{Varadarajan1958}
Varadarajan VS (1958) On the convergence of sample probability distributions.
  Sankhyā: The Indian Journal of Statistics (1933-1960) 19(1/2):23--26

\bibitem[{Villani(2008)}]{Villani2008}
Villani C (2008) Optimal Transport: Old and New. Springer Berlin Heidelberg

\bibitem[{Wainwright(2019)}]{Wainwright2019}
Wainwright MJ (2019) High-Dimensional Statistics: A Non-Asymptotic Viewpoint.
  Cambridge Series in Statistical and Probabilistic Mathematics, Cambridge
  University Press

\bibitem[{von Waldenfels(1968)}]{vonWaldenfels1968}
von Waldenfels W (1968) Characteristische funktionale zufaelliger masse.
  Zeitschrift f{\"u}r Wahrscheinlichkeitstheorie und verwandte Gebiete
  10:279–283

\bibitem[{Wasserman(2006)}]{Wasserman2006}
Wasserman L (2006) All of Nonparametric Statistics (Springer Texts in
  Statistics). Springer-Verlag, Berlin, Heidelberg

\bibitem[{Weed and Bach(2019)}]{WeedBach2019}
Weed J, Bach F (2019) {Sharp asymptotic and finite-sample rates of convergence
  of empirical measures in Wasserstein distance}. Bernoulli 25(4A):2620 -- 2648

\bibitem[{Yurochkin et~al(2019)Yurochkin, Claici, Chien, Mirzazadeh, and
  Solomon}]{Yurochkin2019}
Yurochkin M, Claici S, Chien E, et~al (2019) Hierarchical optimal transport for
  document representation. In: Wallach H, Larochelle H, Beygelzimer A, et~al
  (eds) Advances in Neural Information Processing Systems, vol~32. Curran
  Associates, Inc.

\bibitem[{Zolotarev(1984)}]{Zolotarev1984}
Zolotarev VM (1984) Probability metrics. Theory of Probability and its
  Applications 28:278--302

\end{thebibliography}

\end{document}